\documentclass[twoside]{article}
\usepackage[a4paper]{geometry}
\usepackage[latin1]{inputenc} 
\usepackage[T1]{fontenc} 
\usepackage{RR}
\usepackage{hyperref}
\RRNo{7003}
\RTNo{0703}
\RRdate{November 2012}

\RRauthor{
Jean-Fr\'{e}d\'{e}ric Gerbeau
  \thanks[sfn]{Inria Paris-Rocquencourt, B.P. 105 Domain de Voluceau, 78153 Le Chesnay, France}%
	\thanks[sfnp]{UPMC Univ Paris 06, Laboratoire Jacques-Louis Lions, F-75005, Paris, France}
  \and
Damiano Lombardi
\thanksref{sfn}\thanksref{sfnp}
}

\authorhead{J-F. Gerbeau \& D. Lombardi}
\RRtitle{R\'eduction de mod\`eles bas\'ee sur des paires de Lax approch\'ees}
\RRetitle{Reduced-Order Modeling based on Approximated Lax Pairs}
\titlehead{Reduced-Order Modeling based on Approximated Lax Pairs}
\RRresume{Nous proposons une m\'ethode de r\'eduction de mod\`eles, bas\'ee sur des paires de Lax approch\'ees, pour r\'esoudre des \'equations d'\'evolution non lin\'eaires. Contrairement \`a d'autres m\'ethodes de r\'eduction de mod\`eles, comme la POD, l'espace dans lequel la solution est cherch\'ee \'evolue selon une dynamique reli\'ee au probl\`eme. La m\'ethode est par cons\'equent bien adapt\'ee \`a des probl\`emes comportant des ondes progressives ou des propagations de fronts. Nous montrons des exemples num\'eriques pour les \'equations de KdV et FKPP en dimension un et deux. 
}
\RRabstract{
A reduced-order model algorithm, based on approximations of Lax pairs, is proposed to solve nonlinear evolution partial differential equations. Contrary to other reduced-order methods, like Proper Orthogonal Decomposition, the space where the solution is searched for evolves according to a dynamics specific to the problem. It is therefore well-suited to solving problems with progressive waves or front propagation. Numerical examples are shown for the KdV and FKPP (nonlinear reaction diffusion) equations, in one and two dimensions. }
\RRkeyword{Reduced Order Modeling, Lax Pair, Scattering Transform, KdV, FKPP}
\RRmotcle{R\'eduction de mod\`ele, paires de Lax, transform\'ee de scattering, KdV, FKPP}
\RRprojet{Reo}
\RCParis 

\usepackage{graphicx}
\usepackage{amsmath}
\usepackage{amssymb}
\usepackage{amsthm}
\usepackage{bbold}

\usepackage{color}
\usepackage[normalem]{ulem}
\definecolor{myorange}{rgb}{0.9568,0.4941,0.1961}
\definecolor{myred}{rgb}{0.9098,0.1294,0.2078}
\definecolor{myblue}{rgb}{0.0352,0.4981,0.6509}
\definecolor{mygreen}{rgb}{0.2235,0.6353,0.2588}
\usepackage{marvosym}

\newcommand{\R}{\ensuremath{\mathbb{R}}}

\theoremstyle{plain}
\newtheorem{prop}{Proposition}
\theoremstyle{definition}
\newtheorem*{rem}{Remark}

\begin{document}
\RRNo{8137}
\makeRR   
\section{Introduction}

This work is devoted to a method for solving nonlinear Partial Differential Equations (PDEs) by a reduced-order model (ROM) based on the concept of Lax pairs. 

In his seminal work~\cite{lax-68}, Lax proposed a formalism to integrate a class of nonlinear evolution PDEs. He introduced a pair of linear operators $\mathcal{L}(u)$ and $\mathcal{M}(u)$, where $u$ denotes the solution of the PDE. 
The eigenfunctions of $\mathcal{L}(u)$, are propagated by a linear PDE involving $\mathcal{M}(u)$. For some PDEs,  the eigenvalues of $\mathcal{L}(u)$ have the remarkable property to be constant in time, as soon as  $\mathcal{L}(u)$ and $\mathcal{M}(u)$ satisfy the Lax equation (see \eqref{eq:lax-equation} below). This formalism, which is closely related to the inverse scattering method, can be applied to a wide range of problems arising in many fields of physics, such as Korteweg-de Vries (KdV), Camassa-Holm, Sine-Gordon or nonlinear Schr\"{o}dinger equations. 

The main idea of the present paper is to use the eigenfunctions of $\mathcal{L}(u)$ as a basis to approximate the solutions of the PDE. Contrary to standard ROM techniques, like the Proper Orthogonal Decomposition (POD, see e.g. \cite{Sirovich_1}) or the Reduced Basis Method (RBM, see e.g. \cite{maday-ronquist-02,rozza-huynh-patera-07}), the basis is therefore time dependent. This property makes the method especially well-suited to problems featuring propagation phenomena, which are known to be difficult to tackle with POD.

To define a model that is genuinely ``reduced'', the number of eigenfunctions used to approximate the solution has to be small. This seems to be the case for various problems of practical interest. This has been recently pointed out by Laleg, Cr\'epeau and Sorine who used the eigenfunctions of a Schr\"odinger operator, playing the role of $\mathcal{L}(u)$, 
to provide a parsimonious representation of the arterial blood pressure~\cite{Sorine_1,laleg-08,laleg-crepeau-sorine-12}. Their signal processing technique, called SCSA (for Semi-Classical Signal Analysis), has been the starting point of our work. 

In the huge literature devoted to Lax pairs (see for example \cite{Drazin_1,hickman-hereman-goktas-12} and the reference therein), operator $\mathcal{M}(u)$ is generally used, or searched for, in closed-form. One of the contributions of the present work is to propose an approximation of $\mathcal{M}(u)$ that can be used even when the Lax pair is not explicitely known. In addition, whereas most of the studies consider one dimensional domains and functions rapidly decreasing at infinity or periodic boundary conditions (see \cite{Fokas_1} for a theory on the finite line), our approximation method can easily be used on multi-dimensional bounded domains, with standard boundary conditions. 

The structure of the work is as follows: in Section~\ref{sec:inv-scat-PDE}, the Lax pair is introduced for the KdV equation and the link with the inverse Scattering Transform method is recalled; in Section~\ref{sec:modal-scatter}, numerical approximations of the Lax operators are proposed; Section~\ref{sec:alp-rom} is devoted to our reduced-order model algorithm, that will be called ALP, for Approximated Lax Pair. Some numerical tests are presented in Section~\ref{sec:first-tests} for the Korteweg-de Vries and Fisher-Kolmogorov-Petrovski-Piskunov equations in one and two dimensions. 
\section{The inverse scattering method for nonlinear PDEs}
\label{sec:inv-scat-PDE}


Lax method has strong connections with the inverse scattering method, which can be roughly viewed as the counterpart of the 
Fourier transform for nonlinear problems (see e.g. \cite{Ablowitz_1,Drazin_1}). In this section, we briefly recall how the scattering transform is classically used to solve nonlinear PDEs. This reminder is not strictly necessary to introduce our method, but the reader might find it useful to understand the underlying ideas. We consider the standard example of Korteweg-de Vries equation (KdV):
\begin{equation}
	\label{eq:kdv}
\partial_t u + 6 u\partial_x u + \partial^3_{xxx} u = 0,
\end{equation}
for $(x,t)\in \R\times \R_+$, with the initial condition $u(x,0) = u_0(x), x\in\R$. The Lax pair associated with this equation is defined by
\begin{equation}
\begin{split}
& \mathcal{L}(u) v = -\partial^2_{xx} v -u v, \\
& \mathcal{M}(u) v = 4 \partial^3_{xxx} v +6 u\partial_x v + 3 v\partial_x u. \\
\end{split}
\label{KdV_Operators}
\end{equation}
These operators satisfies the Lax equation
\begin{equation}
	\label{eq:lax-equation}
	\partial_t \mathcal{L}(u) + \mathcal{L}(u)\mathcal{M}(u) - \mathcal{M}(u)\mathcal{L}(u) = 0,
\end{equation}
if and only if $u(x,t)$ is solution to the KdV equation.
Consider the spectral problems, parametrized by $t$,
\[
\mathcal{L}(u(x,t))\phi(x,t) = \lambda(t)\phi(x,t).
\] 
It can easily be deduced from the Lax equation that 
\begin{equation}
	\label{eq:isospectral}
	\partial_t \lambda = 0.
\end{equation}
When this property is satisfied, the problem is said to be \emph{isospectral}. 
The Lax equation describes how the operator $\mathcal{L}(u)$ evolves in time. Its eigenfunctions can be used to reconstruct the solution $u$ at every time. 


The solution of the KdV equation by inverse scattering can be decomposed in three steps. 

The first step is to solve the spectral problem associated with the initial data:
\[
\mathcal{L}(u_0)\phi = \lambda \phi
\]
which is the linear Schr\"odinger equation with potential $-u_0$. In general, this problem, known as the \emph{Direct Scattering Transform}, has a continuous spectrum of positive values and a discrete spectrum of negative eigenvalues. For the sake of simplicity, we assume that $u_0$ is such that it can be reconstructed with the sole discrete part of the spectrum  (such a $u_0$ is known as a \emph{reflectionless potential}). We also assume that $u_0$ is in the Schwartz space $\cal{S}(\R)$, i.e. decays sufficiently rapidly at infinity, which implies that the number of negative eigenvalues is finite. For $m=1,\cdots,N_-$, we denote by $\phi_m$ the eigenfunctions normalized in $L^2(\R)$, by $\lambda_m$ the eigenvalues, we define $\kappa_m = \sqrt{-\lambda_m}$ and we rank the eigenvalues such that $\kappa_1 > \kappa_2 > \cdots > \kappa_{N_-}$. The Deift-Trubowitz formula then gives (see \cite{Ablowitz_1}):
\begin{equation}
	\label{eq:u0-deift-trubo}
u_0(x) = \sum_{m=1}^{N_-}\kappa_m \phi_m^2(x,0).
\end{equation}
When $x$ goes to infinity, the eigenfunctions $\phi_m(x,0)$ can be proved to be equivalent to $c_m(0)\exp(- \kappa_m x)$. The quantities $(\kappa_m(0), c_m(0))$ are known as the \emph{scattering data}. Note that in general, the scattering data also contain the so-called \emph{reflexion coefficient}, which is not present here because of the assumption made on $u_0$. This missing part is important for scattering problems, but not for the present work. 


The second step consists in propagating the scattering data in time. This is trivial for the $\kappa_m$, since they are already known to be independent of $t$ (see \eqref{eq:isospectral}). For $c_m$, it can be shown, using the fact that the solution is in the Schwartz space, that 
\begin{equation}
	\label{eq:edo-ck}
	\frac{dc_m}{dt}(t) = 4 \kappa_m^3 c_m(t),	
\end{equation}
which readily gives $c_m(t)$. 

The third step is the \emph{Inverse Scattering Transform}, i.e. compute $u(x,t)$ from the scattering data $(\kappa_m,c_m(t))$. In general, this can be done by solving an integral equation, known as the Gelfand-Levitan-Marchenko equation. When $u_0$ is reflectionless, \emph{i.e.} given by~\eqref{eq:u0-deift-trubo}, a closed-form expression of $u(x,t)$ using the scattering data is known. More details can be found e.g. in \cite{Ablowitz_1,Drazin_1}.
 
The method that has just been described heavily relies on the specificity of the KdV equation set in the whole space $\R$, with an initial condition rapidly decreasing at infinity. In the remainder of this article, we propose a method directly based on an approximation of the Lax pair. It formally follows the same steps as the inverse scattering method, but can be used in a more general setting. 

\section{Preliminary results}
\label{sec:modal-scatter}

We consider an evolution PDE in a bounded subset $\Omega$ of $\mathbb{R}^d$, $d\geq 1$, of the form:
\begin{equation}
	\label{DiffProb_1}
	\partial_{t} u = F(u),
\end{equation}
where $F(u)$ is an expression involving $u$ and its derivatives with respect to ${x}_1,\dots,{x}_d$. The problem is completed with an initial condition 
\begin{equation}
	\label{DiffProb_IC}
	u({x},0) = u_0({x}), \mbox{ for } {x}=({x}_1,\dots,{x}_d) \in \Omega.
\end{equation}
For the sake of simplicity in this Section, $u$ is assumed to vanish on the boundary $\partial\Omega$. Other boundary conditions will be considered in the numerical tests.

The solution to this problem is searched for in an Hilbert space $V$ (of functions vanishing on $\partial\Omega$) and approximated in $V_h$, a finite dimensional subspace of $V$, for example obtained by the finite element method. Let $(v_{j,h})_{j=1..N_{h}}$ denotes a basis of $V_h$ and 
$\langle\cdot,\cdot\rangle$ the $L^2(\Omega)$ scalar product. 

We follow the three steps presented in the previous section to solve the KdV equation by inverse scattering. For each step, we propose a numerical approximation that will be used for the reduced order model integration. 

\subsection{Semi Classical Signal Analysis}
\label{sec:scsa}

Consider a nonnegative signal\footnote{If $u({x})$ is not nonnegative, it is replaced by $u({x}) - \min_{x\in\Omega} u_0(x)$}$u({x})$, for ${x}=({x}_1,\dots,{x}_d)\in\Omega$, a real number $\chi>0$, and the linear Schr\"odinger operator 
\begin{equation}
	\label{eq:def-Lu0-chi}
\mathcal{L}_{\chi}(u) \phi = - \Delta \phi - \chi u \phi,
\end{equation}
where $\Delta$ denotes the Laplacian in $d$ dimensions. The approximated scattering transform we consider is called SCSA for Semi Classical Signal Analysis. It has been proposed in \cite{laleg-08}, analyzed in \cite{laleg-crepeau-sorine-12}, and successfully used for different applications to signal analysis in hemodynamics~\cite{laleg-medigue-sorine-07,laleg-medigue-van-de-louw-10}. It partially relies on the results proved in \cite{Levermore_1,Levermore_2} and it consists of only keeping the modes corresponding to the negative eigenvalues $(\lambda_n)_{n=1... N_-}$ to approximate $u$ by the Deift-Trubowitz formula:
\begin{equation}
	\label{eq:u-deift-trubo}
	\tilde u({x}) = \chi^{-1}\sum_{m=1}^{N_-} \kappa_m \phi_m^2,
\end{equation}
with $\kappa_m = \sqrt{-\lambda_m}$. The parameter $\chi>0$ is chosen in order to reach the desired accuracy. For large values of $\chi>0$, the representation is more accurate, but also more expensive since the number of negative eigenvalues is larger. 

Numerically, we search for $\phi_{m,h} \in V_h$ and $\lambda_{m,h}$ such that
\begin{equation}
	\label{eq:scsa-h}
	\langle \nabla\phi_{m,h},\nabla v_{i,h}   \rangle - \chi \langle u \phi_{m,h}, v_{i,h} \rangle = \lambda_{m,h} \langle \phi_{m,h},v_{i,h}\rangle, \mbox{ for } i=1,\dots,N_{h},	
\end{equation}
and $u$ is approximated by 
\[
\tilde u_h({x}) = \chi^{-1}\sum_{m=1}^{N_-} \kappa_{m,h} \phi_{m,h}^2({x}),
\]
where $\chi$ is chosen such that 
\[
\| u - \tilde u_h \|_{L^2(\Omega)} \leq \epsilon_0,
\]
where $\epsilon_0$ is a prescribed tolerance. 

\begin{rem}
Note that problem~\eqref{eq:scsa-h} implies:
\[
\langle u \phi_{m,h}, \phi_{p,h} \rangle  - \frac{1}{\chi}	\langle \nabla\phi_{m,h},\nabla \phi_{p,h}   \rangle  
+ \frac{\lambda_{m,h}}{\chi}\langle \phi_{m,h},\phi_{p,h}\rangle
=  0, \mbox{ for } p=1,\dots,N_{-},	
\]
which is the Euler-Lagrange equation of the maximization of the $L^2$ projection of $u$ onto the space spanned by the $\phi^2_{m,h}$, up to a regularization term, and under the constraint of a unitary $L^2$ norm.
\end{rem}

\subsection{Modal approximation of the Lax propagation operator}
\label{sec:modalApp}

For an arbitrary PDE \eqref{DiffProb_1}, the construction of $\mathcal{L}(u)$ and $\mathcal{M}(u)$ satisfying 
\begin{equation}
\begin{split}
& \mathcal{L}(u) \phi = \lambda \phi, \\
& \partial_{t} \phi = \mathcal{M}(u) \phi,\\
\end{split}
\label{ScattProb_1}
\end{equation}
and the Lax equation \eqref{eq:lax-equation} is not obvious. In addition, the isospectral property~\eqref{eq:isospectral}, which was instrumental in solving the KdV equation, is in general not satisfied.

In this work, we choose $\mathcal{L}(u)$ as the linear Schr\"odinger operator associated to the potential $-\chi u$, as in \eqref{eq:def-Lu0-chi}. The following proposition shows that it is possible to compute an approximation of $\mathcal{M}(u)$ in the space defined by the eigenfunctions of $\mathcal{L}(u)$ and to derive an evolution equation satisfied by the eigenvalues of $\mathcal{L}(u)$, even when the operator $\mathcal{M}(u)$ is not known in a closed-form. 

From now on, the eigenfunctions of $\mathcal{L}(u({t}))$ will be denoted by $(\phi_m({t}))_{m=1 \cdots N_-}$ when they are used to approximate the solution $u({t})$, and by $(\psi_m({t}))_{m=1 \cdots N_M}$ when they are used to approximate the operator $\mathcal{M}(u({t}))$ or the evolution equation satisfied by the eigenvalues. Functions $\phi_m$ and $\psi_m$ are therefore the same, but the numbers $N_-$ and $N_M$ will be different in general: while the solution will be approximated only on modes corresponding to negative eigenvalues, the approximation of the operators will also necessitate some modes associated with positive eigenvalues. Thus $N_M$ will be larger than $N_-$ in practice.

\begin{prop}
	\label{prop:m-lambda}
	Let $u$ be a solution of equation \eqref{DiffProb_1}. Let $\mathcal{L}_\chi(u)$ be defined by	
	\begin{equation}
		\label{eq:def-Lu-chi}
	\mathcal{L}_\chi(u) \psi = - \Delta \psi - \chi u \psi	
	\end{equation}
    where $\chi$ is a given positive real number. 

Let $N_M\in\mathbb{N}^\ast$. For $m\in \{1,\dots,N_M\}$, let $\lambda_m({t})$ be an eigenvalue of $\mathcal{L}_\chi(u({x},{t}))$, and $\psi_m({x},{t})$ an associated eigenfunction, normalized in $L^2(\Omega)$. Assume there exists an operator $\mathcal{M}(u)$ such that 
	\begin{equation}
	\partial_{{t}} \psi_m = \mathcal{M}(u) \psi_m.
	\label{eq_dtPhi}
	\end{equation}
Then the evolution of $\lambda_m$ is governed by
\begin{equation}
	\partial_{{t}} \lambda_m = - \chi \langle  F(u)\psi_m, \psi_m\rangle,
\label{lambda_eq}
\end{equation}
and the evolution of $\psi_m$ satifies, for $p\in\{1,\dots,N_M\}$,
\begin{equation}
	\label{eq:mode_eq}
	\langle \partial_{t} \psi_m, \psi_p \rangle  = M_{mp}(u),
\end{equation}
with
\begin{equation}
	\label{M_eq}
	\left\{
	\begin{array}{rcl}
M_{mp}(u) &=& \displaystyle{\frac{\chi}{\lambda_p - \lambda_m}} \langle  F(u)\psi_m, \psi_p\rangle, ~~\mbox{ if } p\neq m \mbox{ and } \lambda_p \neq \lambda_m,\\
M_{mp}(u) &=& 0, ~~ \mbox{ if } p=m  \mbox{ or } \lambda_p = \lambda_m.
\end{array}
\right.
\end{equation}
We will denote by $M(u) \in \R^{N_M\times N_M}$ the skew-symmetric matrix whose entries are defined by $M_{mp}(u)$.
\end{prop}

\textbf{Proof.}
Differentiating with respect to $t$ the equation satisfied by the $m$-th mode
\[
\mathcal{L}_\chi(u(x,t)) \psi_m(x,t) = \lambda_m(t) \psi_m(x,t),
\]
we get
\begin{equation}
\left(\mathcal{L}_\chi(u) - \lambda_m \mathcal{I}\right)\partial_{{t}} \psi_m = \partial_{{t}} \lambda_m \psi_m + \chi F(u)\psi_m.
\end{equation}

The scalar product is taken with a generic $\psi_p$, leading to:
\begin{equation}
\langle  \left(\mathcal{L}_\chi(u) - \lambda_m \mathcal{I}\right)\partial_{{t}} \psi_m,\psi_p\rangle  = \partial_{{t}} \lambda_m \langle  \psi_m,\psi_p\rangle  + \chi \langle  F(u)\psi_m, \psi_p\rangle ,
\end{equation}
Using the self-adjointness of the operator and the orthonormality of the eigenfunctions, the following problem is obtained:
\begin{equation}
(\lambda_p-\lambda_m)\langle  \partial_{{t}} \psi_m, \psi_p\rangle  = \partial_{{t}} \lambda_m \delta_{mp} + \chi\langle  F(u)\psi_m, \psi_p\rangle .
\label{eq_Lax_1}
\end{equation}
Taking $p=m$, this proves \eqref{lambda_eq}. In addition, the $L^2$ norm of $\psi_m$ being 1, $\langle \partial_{t} \psi_m, \psi_p \rangle = 0$, i.e. \eqref{M_eq}$_2$.

If $p\neq m$, but $\lambda_p = \lambda_m$ (multiple eigenvalues), we arbitrarilly set $M_{mp}(u) = 0$.

For $\lambda_p\neq \lambda_m$, we deduce from~\eqref{eq_Lax_1}:
\begin{equation}
\langle  \partial_{{t}} \psi_m, \psi_p\rangle  =  \frac{\chi }{\lambda_p - \lambda_m} \langle  F(u)\psi_m, \psi_p\rangle.
\label{dtMode}
\end{equation}
Combining \eqref{eq_dtPhi} with \eqref{dtMode}, equation \eqref{M_eq} is obtained, which completes the proof.
$\diamondsuit$

Equation \eqref{M_eq} gives an approximation of the operator $\mathcal{M}(u)$ on the basis defined by the modes at time ${t}$. This representation is convenient from a computational standpoint since it can easily be obtained from the expression $F(u)$ defining the PDE~\eqref{DiffProb_1}, without any \emph{a priori} knowledge of $\mathcal{M}(u)$. With this approximation of $\mathcal{M}(u)$, the evolution of the modes can be computed according to the nonlinear dynamics of the system. This is an important difference with standard reduced-order methods, like POD, where the modes are fixed once for all.

To set up a reduced order integration method, only a small number $N_M$ of modes will be retained. This number has to be chosen in order to represent the dynamics in a satisfactory way. A possible indicator of the quality of the approximation is given by the following quantity
\begin{equation}
\label{def-norm-dt-psi}
e(\psi_m({t}), N_M)  = \sum_{n=1}^{N_M} (M_{mn}(u({t})))^2,
\end{equation} 
which is an approximation of the $L^2$ norm of the time derivative of $\psi_m$:
\[
\int_{\Omega} \left(\partial_{{t}} \psi_m\right)^2\ d\Omega \approx \sum_{n,l=1}^{N_M} M_{ml}(u)M_{mn}(u) \langle  \psi_n,\psi_l\rangle  = \sum_{n=1}^{N_M} M_{mn}(u)^2 = e(\psi_m, N_M).
\]

By summing up over the modes, the Frobenius norm of the representation of the evolution operator is recovered:
\begin{equation}
\label{frobDef}
\sum_m^{N_M} e(\psi_m({t}), N_M)  = \sum_{m,n=1}^{N_M} (M_{mn}(u({t})))^2 = F_{N_M}^2.
\end{equation} 

This norm may be used as an error indicator for the dynamics recovery and it was investigated by the numerical experiments presented in Section~\ref{sec:first-tests}.

\subsection{Approximated Inverse Scattering Transform}
\label{sec:approx-ist}

Proposition~\ref{prop:m-lambda} gives an approximated way to propagate the eigenmodes and the eigenvalues associated to a Lax pair. From a given set of eigenmodes and eigenvalues, we will have to reconstruct the solution. We therefore need an approximated counterpart of the inverse scattering step presented in Section~\ref{sec:inv-scat-PDE}.   

Assume that a family $(\phi_m)_{m=1..N_-}$ of eigenfunctions of $\mathcal{L}_\chi(u)$, and eigenvalues $(\lambda_m)_{m=1..N_-}$ are known. We propose to approximate $u$ as
\[
\tilde u = \sum_{k=1}^{N_-} \alpha_k \phi_k^2.
\]
Inserting this expression in
\[
\langle \mathcal{L}_\chi(\tilde u) \phi_m,\phi_p\rangle = \lambda_m \langle \phi_m, \phi_p \rangle,
\]
and using that $\langle \phi_m, \phi_p \rangle = \delta_{mp}$, we obtain:
\begin{equation}
	\label{eq:approx-ist}
\sum_{k=1}^{N_{-}} \alpha_k \langle  \phi_k^2,\phi_m^2 \rangle  =- \frac{1}{\chi} \left( \lambda_m + \langle  \Delta\phi_m, \phi_m\rangle \right).
\end{equation}
This is a linear system for $\alpha_k$, whose resolution is costless when a small number of modes is considered.


\section{Reduced-Order Modeling based on Approximated Lax Pairs (ALP)}
\label{sec:alp-rom}

The eigenvalues and the eigenfunctions of $\mathcal{L}_\chi(u(\cdot,{t}))$ evolve as ${t}$ changes. We propose in this section numerical schemes to approximate these evolutions. 

\subsection{Approximated propagation of the eigenvalues}

According to Proposition~\ref{prop:m-lambda}, the evolution of the eigenvalues is governed by~\eqref{lambda_eq}. Numerically, we consider a simple explicit Euler scheme:

Assume that $u^n,\lambda^n$ and $(\psi^n_p)_{p=1,\dots,N_M}$, are known at step $n$, compute $\lambda_p^{n+1}$ as:
\begin{equation}
\lambda_p^{n+1} = \lambda_p^{n} - \delta{t} \langle  F(u^n)\psi_p^n,\psi_p^n\rangle ,
\label{eigenvalUpdate}
\end{equation}
where $\delta{t}$ denotes the time step.

\subsection{Approximated propagation of the eigenfunctions}

To propagate the modes $\psi_p$, it is important to devise a procedure that preserves their orthonormality. We first note that equation~\eqref{eq:mode_eq} yields:
\begin{equation}
\label{eq:psi-with-residual}
\frac{\partial\psi_m}{\partial{t}} = \sum_{p=1}^{N_M} M_{mp}(u) \psi_p\ + r_m,
\end{equation}
where $r_m({t})\in\left[\mbox{span}(\psi_1({t}),\dots,\psi_{N_M}({t}))\right]^\perp$. 

Assume that $u^n$ and $(\psi^n_p)_{p=1,\dots,N_M}$ are known and that $\langle \psi^n_m, \psi^n_p\rangle = \delta_{mp}$ for $m,p=1,\dots,N_M$ . Consider the following problem, obtained from~\eqref{eq:psi-with-residual} by neglecting the residual $r_m$ and linearizing the first term in the right-hand side:
\begin{equation}
	\label{eq:psi-tilde}
\left\{
\begin{split}
	\frac{\partial\tilde\psi_m}{\partial{t}} &= \sum_{p=1}^{N_M} M_{mp}(u^n) \tilde\psi_p,\\
	\tilde\psi_m(0) &= \psi_m^n,
\end{split}
\right.
\end{equation}
whose solution is 
\begin{equation}
	\label{eq:sol-exp}
	\tilde\psi_m({t}) = \sum_{p=1}^{N_M} \left[ \exp(M(u^n) {t}) \right]_{mp} \psi_p^{n}.
\end{equation}
Note that
\begin{equation}
\langle  \tilde\psi_m({t}),\tilde\psi_p({t})\rangle  =  \left[ \exp(M(u^n){t})\exp(M(u^n)^{T}{t})\right]_{mp} = \delta_{mp},
\end{equation}
since $M(u^n)$ is skew-symmetric. Thus, by construction, this procedure preserves the orthonormality. Ideally, $\psi^{n+1}_m$ could be defined as $\tilde\psi(\delta{t})$. But to avoid the expensive computation of the exponential of a matrix, we suggest to approximate $\exp(M(u^n)\delta{t})$ by 
$T(u^n) = I + M(u^n)\delta{t} + \frac{1}{2}M(u^n)^2\delta{t}^2$. Thus, $\psi^{n+1}_m$ is eventually defined by:
\begin{equation}
	\label{eq:sol-exp-approx}
	 \psi^{n+1}_m = \sum_{p=1}^{N_M} \left[ I + M(u^n)\delta{t} + \frac{\delta{t^2}}{2}M(u^n)^2 \right]_{mp} \psi_p^{n}.
\end{equation}

With this approximation, the orthonormality of $(\psi^{n+1}_m)_{m=1,\dots,N_M}$ is only perturbed with an error of the order of $\delta{t}^4$.

\subsection{The ALP algorithm}

We now have all the tools necessary to the definition of a reduced order model based on the approximated Lax pairs.

\textbf{Initialization:}
Let $u_0$ be the initial condition and let  $\epsilon_0>0$ be a prescribed tolerance. Compute a set of modes $(\psi^0_{m,h})_{m=1\dots N_M}$ and eigenvalues $(\lambda^0_{m,h})_{m=1\dots N_M}$ by solving
\[
\langle \nabla\psi^0_{m,h},\nabla v_{i,h}   \rangle - \chi \langle u \psi^0_{m,h}, v_{i,h} \rangle = \lambda^0_{m,h} \langle \psi^0_{m,h},v_{i,h}\rangle, \mbox{ for } i=1,\dots,N_{h},	
\]
The eigenvalues $\lambda^0_{m,h}$ are ranked in the ascending order. Let $N_-\leq N_M$ be the number of negative eigenvalues. Throughout the algorithm, we define $\phi^n_{m,h} = \psi^n_{m,h}$, for $m\in\{1,\dots,N_-\}$. The initial condition is approximated by the SCSA method (see Section~\ref{sec:scsa}): 
\[
\tilde u^0_h({x}) = \chi^{-1}\sum_{m=1}^{N_-} \kappa_{m,h} \phi_{m,h}^2({x}),
\]
where $\chi$ is chosen such that $\| u_0 - \tilde u^0_h \|_{L^2(\Omega)} \leq \epsilon_0$.

\textbf{Time evolution:} For $n\geq 0$, assume that $u^{n}_h, (\psi_{m,h}^{n})_{m=1\cdots N_M}$ and $\lambda_k^{n}$ are known at time $t^n$.
\begin{enumerate}
\item Compute the matrix $M(u^{n}_h)$ approximating the operator $\mathcal{M}(u(t^n))$ with~(\ref{M_eq}): 
\[
M_{mp}(u^{n}_h) = \frac{\chi}{\lambda^{n}_p - \lambda^n_m} \langle  F(u^n_h)\psi^n_{m,h}, \psi^n_{p,h}\rangle, ~~\mbox{ for } m,p = 1,\dots,N_M.
\]
\item Compute the eigenvalues $\lambda_{m,h}^{n+1}$ with (\ref{eigenvalUpdate}): 
\[
\lambda_{m,h}^{n+1} = \lambda_{m,h}^{n} - \delta{t} \langle  F(u^n_h)\psi_{m,h}^n,\psi_{m,h}^n\rangle ,
\]
\item Compute the eigenfunctions $\psi_m^{n+1}$ with (\ref{eq:sol-exp-approx}):
\[
\psi^{n+1}_{m,h} = \sum_{p=1}^{N_M} \left[ I + \delta{t}M(u^n_h) + \frac{\delta{t^2}}{2}M(u^n_h)^2 \right]_{mp} \psi_{p,h}^{n}.
\]
\item Solve system~\eqref{eq:approx-ist} for $\alpha^{n+1}_p$:
\[
\sum_{p=1}^{N_{-}} \alpha^{n+1}_p \langle  \phi_{p,h}^2,\phi_{m,h}^2 \rangle  =- \frac{1}{\chi} \left( \lambda^{n+1}_p + \langle  \Delta\phi^{n+1}_{m,h}, \phi^{n+1}_{m,h}\rangle \right).
\]
\item Compute $u^{n+1}_h$ with
\[
u^{n+1}_h = \sum_{p=1}^{N_-} \alpha_p^{n+1} \left(\phi^{n+1}_{p,h}\right)^2.
\]
\end{enumerate}

\begin{rem}
	\label{rem:variable-num-mode}
For the sake of simplicity, $N_-$ and $N_M$ are fixed in this algorithm. They of course could be adapted along the resolution acording to various criteria. We will not investigate this possibility in the present study, except in the following case: if one eigenvalue which was initially positive becomes negative, then the corresponding eigenmode can be simply added to the set used to approximate the solution, and the value of $N_-$ incremented accordingly. This will be used in the test case presented in Section~\ref{sec:fkpp-1d}.
\end{rem}

\section{Numerical Experiments}
\label{sec:first-tests} 

In this section, some numerical experiments are presented. The aim is to assess the proposed technique and to highlight some differences with the POD.

Two test cases deal with the Korteweg-de Vries (KdV) equation. As recalled in Section~\ref{sec:inv-scat-PDE}, this equation is a remarkable example of an integrable system, which can be exactly solved by inverse scattering (see \cite{Drazin_1,Ablowitz_1}). This problem is therefore a good candidate to assess our numerical approach in a case when an analytical solution is possible. 

Then, the Fisher-Kolmogorov-Petrovski-Piskunov (FKPP) equation is considered in one and two dimensions. This equation, which arises in many applications, features front propagation. Contrary to KdV, it is not isospectral, and no closed-form expression of the Lax pair is known. So, it is an interesting problem to test our method when no analytical approach is available.

\subsection{Korteweg-de Vries equation}

We consider the one-dimensional KdV equation \eqref{eq:kdv}, with an initial condition $u_0$ that will be either a one-soliton or a three-soliton. The computational domain is bounded for practical reason, but the time range of the computations and the length of the domain are chosen such that the boundaries have no effect on the solution. 

\subsubsection{One-soliton propagation}
\begin{figure}
\centerline{\hbox{\begin{tabular}{cc}
\includegraphics[height=5.0cm,width=6.5cm]{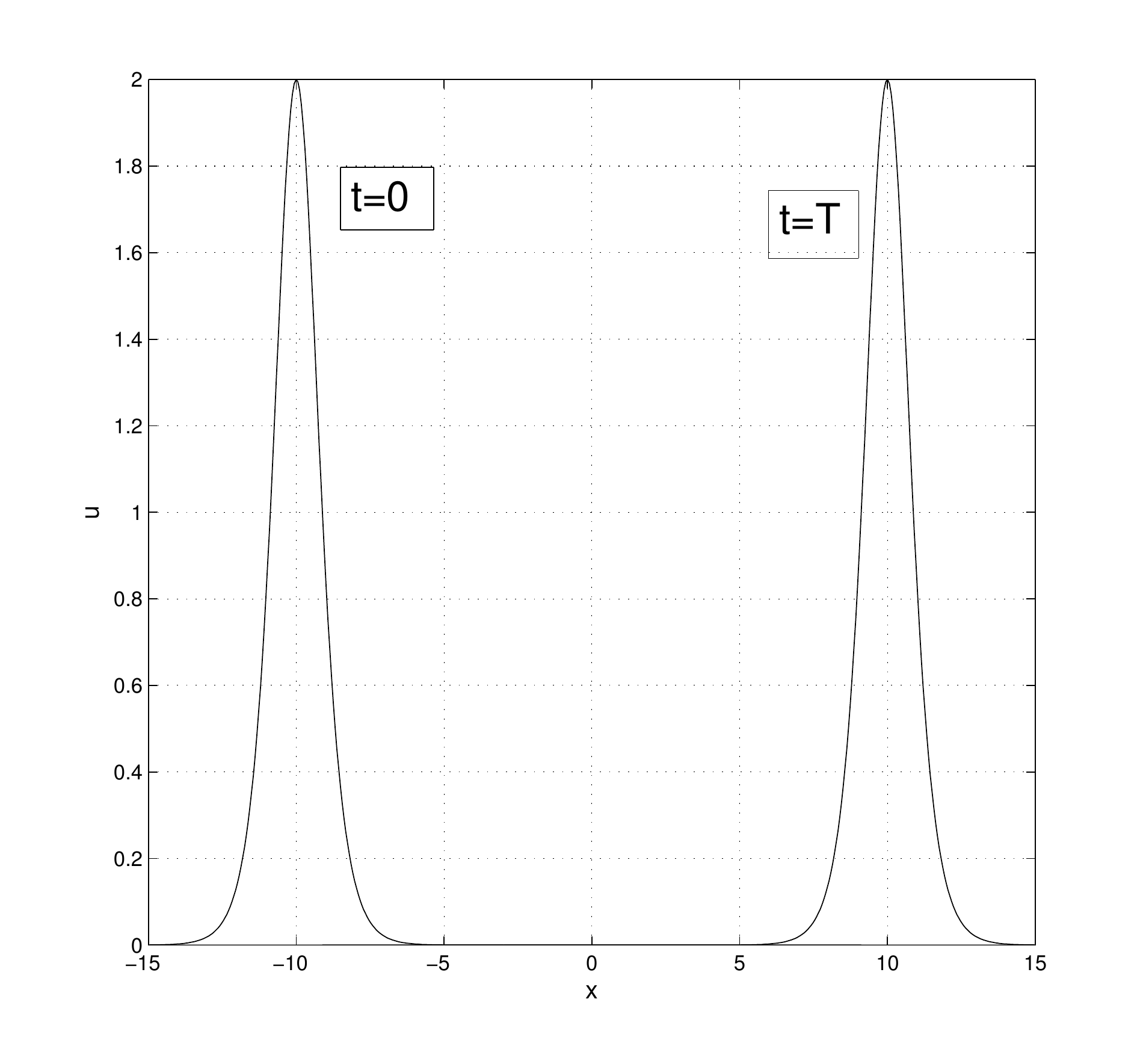} &
\includegraphics[height=5.0cm,width=6.5cm]{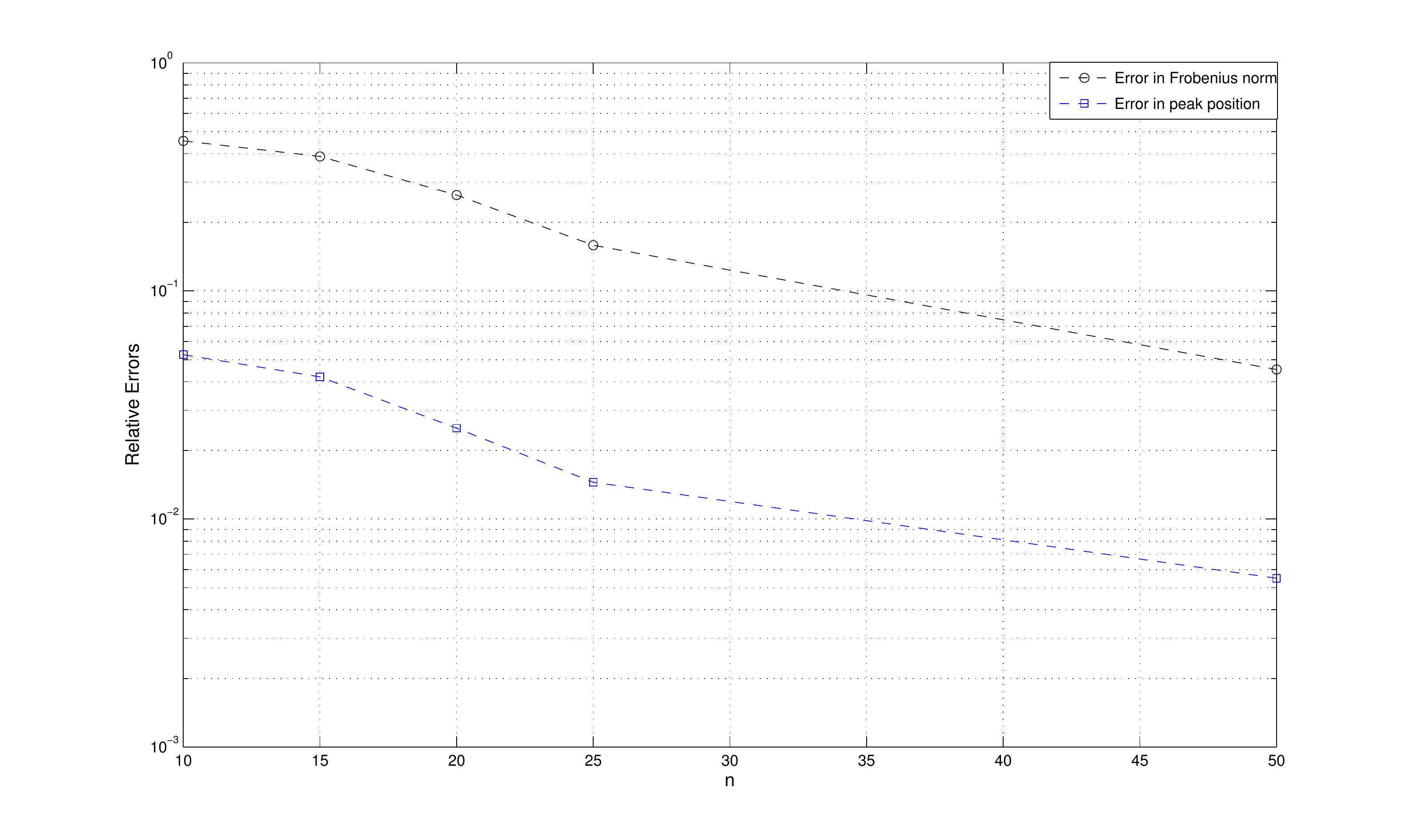}\\
\hspace{0.3cm} (a) & \hspace{0.3cm}(b) 
\end{tabular}}}
\caption{a) Exact one-soliton solution, see Eq.(\ref{KdV_OneSolExact}), at initial and final time (T=5.0) b) Relative error for the Frobenius norm of the operator $M$ and for the peak position, at $t=T$, in semi-logarithmic plot.}
\label{KdV_F1}
\end{figure}

The reference solution shown in Fig.\ref{KdV_F1}.a) at initial and final time ($T=5.0$) is a one-soliton propagation, namely:
\begin{equation}
u(x,t) = \frac{\beta}{2} sech^2\left(\frac{\beta^{1/2}}{2}(x-\beta t) \right),
\label{KdV_OneSolExact}
\end{equation}
with $\beta = 4$.

The modes are extracted by using the initial condition only $u_0=u(x,0)$, setting $\chi = 1$. The Schr\"{o}dinger spectral problem was discretized in a space of $N_{h}=500$ piecewise linear functions.

As $u_0$ is the initial datum of the one-soliton propagation and $\chi=1$ provides the analytical expression for $\mathcal{L}(u)$ in the case of the KdV equation (see Eq.(\ref{KdV_Operators})), only one eigenvalue belongs to the discrete spectrum and the corresponding mode squared is exactly $u_0$, up to discretization errors ($10^{-4}$ in $L^2$ norm for the present case). Hence, only one mode is retained in order to represent the solution ($N_-=1$). 

\begin{figure}
\centerline{\hbox{\begin{tabular}{cc}
\includegraphics[height=5.0cm,width=6.5cm]{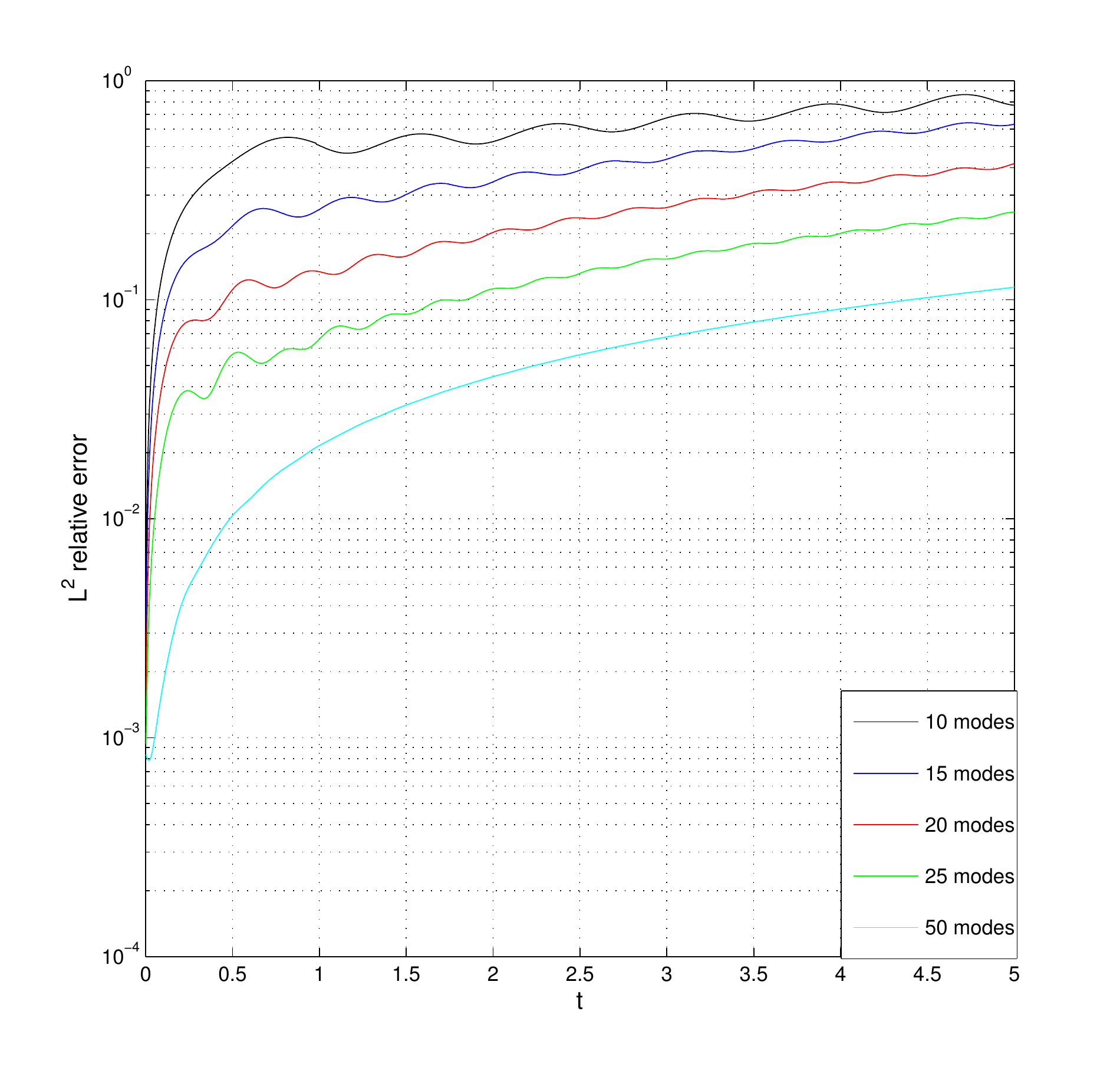} &
\includegraphics[height=5.0cm,width=6.5cm]{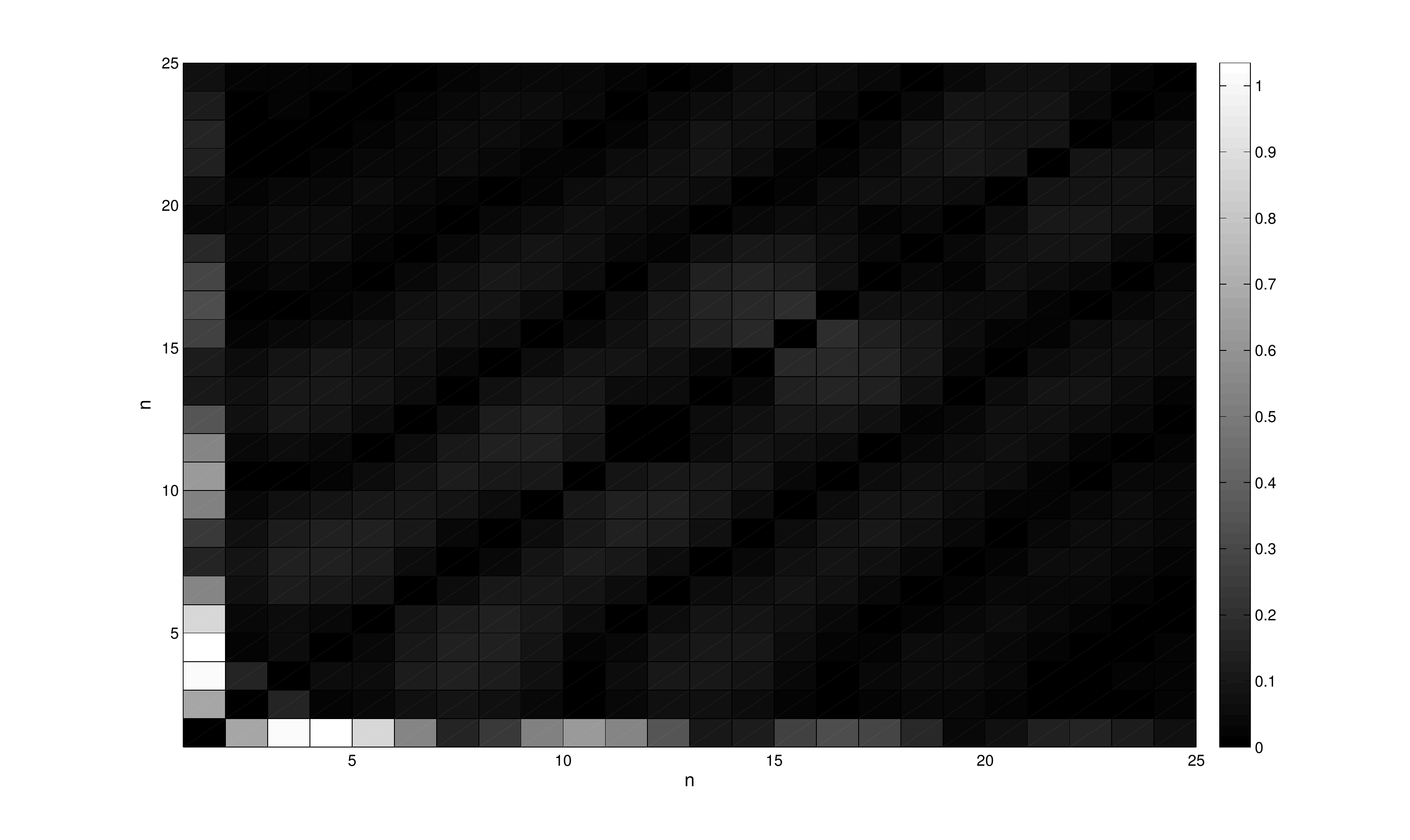}\\
\hspace{0.3cm} (a) & \hspace{0.3cm}(b) 
\end{tabular}}}
\caption{a) Errors in $L^2$ norm while varying the modes number in the approximation of the operator $\mathcal{M}(u)$, b) representation of $\mathcal{M}(u)$, in absolute value.}
\label{KdV_F3}
\end{figure}

The influence of the number $N_M$ of modes on the dynamics approximation is investigated by computing the error in $L^2$ norm of $u$ as a function of time. The results are shown in Fig.\ref{KdV_F3}.a), in a semi-logarithmic plot. As expected, the error globally decreases when the number of modes is increased and it tends to increase in time, due to the accumulation of the errors during the integration. 

A meaningful and synthetic criterion of the dynamics recovery quality is the relative error in the peak position at final time, that is, the difference, in modulus, between the position of the peak of the exact and the reconstructed solution, normalized by the distance traveled by the wave, which is $20$ in this case, see Fig. \ref{KdV_F1}.a). 

It is enlightening to compare the error in the peak position and an indicator of the approximation error of operator $\mathcal{M}(u)$. To estimate this indicator, we compute the relative difference between a reference value and the Frobenius norms (see Section~\ref{sec:modalApp}) of matrices $M(u)$ for $N_{M} = 10,15,20,25,50$, computed at final time $t=T$:
\begin{equation}
	\label{eq:error-frob}
\mathcal{E}_{F}(N_M,N_\infty):=\frac{F_{N_\infty} - F_{N_M} }{F_{N_\infty}},
\end{equation}
where $F_{N} = (\sum_{m,p=1}^{N} M_{mp}^2)^{1/2}$ and where the reference value $N_\infty=100$ (adding more modes changes its value less than $0.5\%$). In Fig. \ref{KdV_F1}.b),  the error in the peak position and $\mathcal{E}_{F}(N_M, N_\infty)$ are shown as a function of $N_M$, in a semi-logarithmic plot. The trend of the curves is practically the same, suggesting that $\mathcal{E}_{F}(N_M,N_\infty)$ might be used as an indicator to assess the quality the approximation.

As said above, the KdV equation defines an isospectral flow, that is, the eigenvalues of the Schr\"{o}dinger operator associated with the solution of the problem do not vary in time. An interesting validation test is to see how far the reduced order approach satisfies the isospectrality property. The computation of the right-hand side of Eq.(\ref{lambda_eq}) provides the right result, for all the number of modes used. When $\delta{t} = 2.5 \ 10^{-3}$ is used for the time discretization, the difference between the right value of the first (negative) eigenvalue (namely $\lambda = -1$) and the computed one is $4.5\ 10^{-5}$ at final time $t=T$.

Let us recall that only one mode is used in order to represent the quantity $u$, while a varying number of modes is adopted to study the dynamics.  
\begin{figure}
\centerline{\hbox{\begin{tabular}{cc}
\includegraphics[height=5.0cm,width=6.5cm]{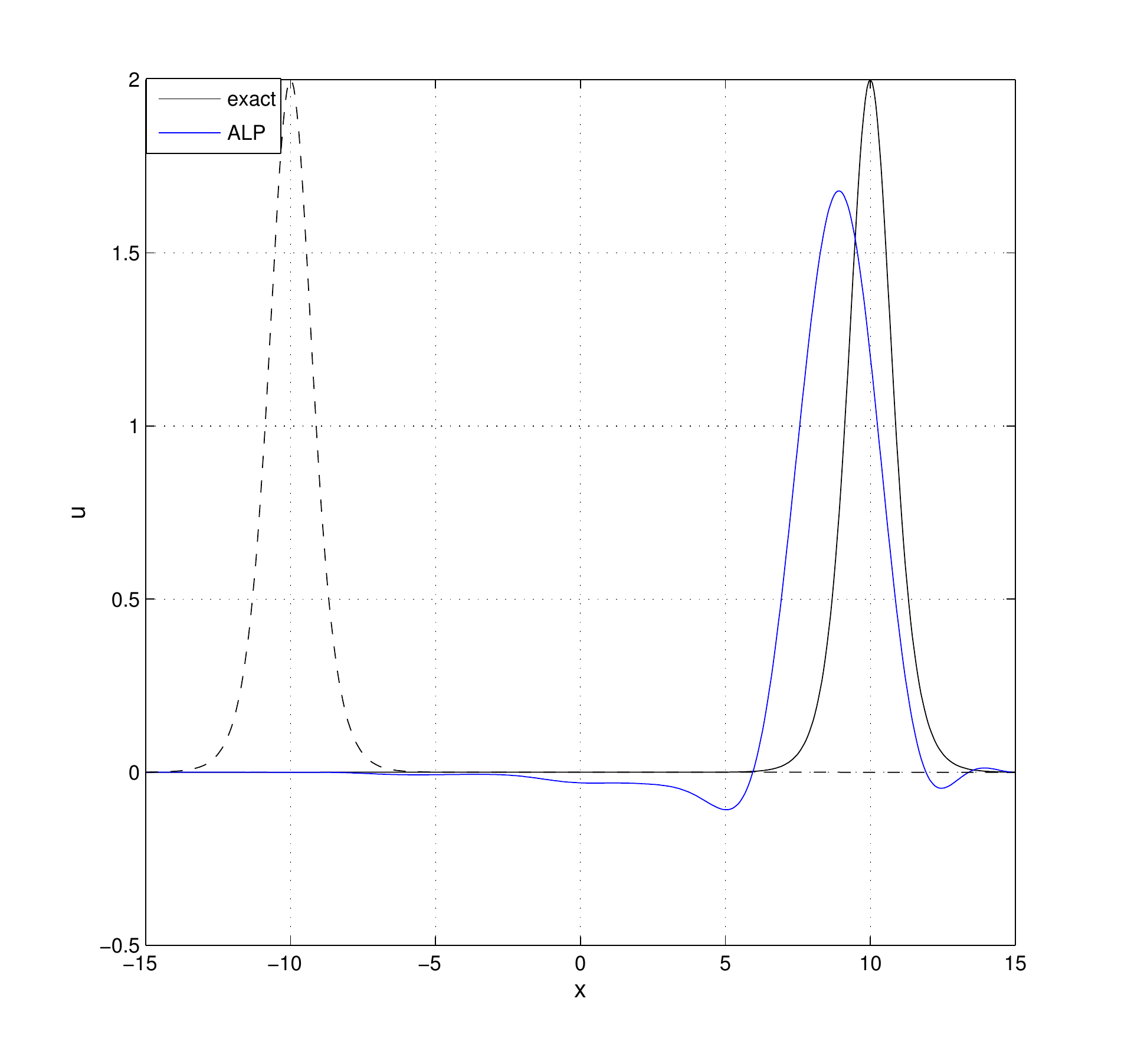} &
\includegraphics[height=5.0cm,width=6.5cm]{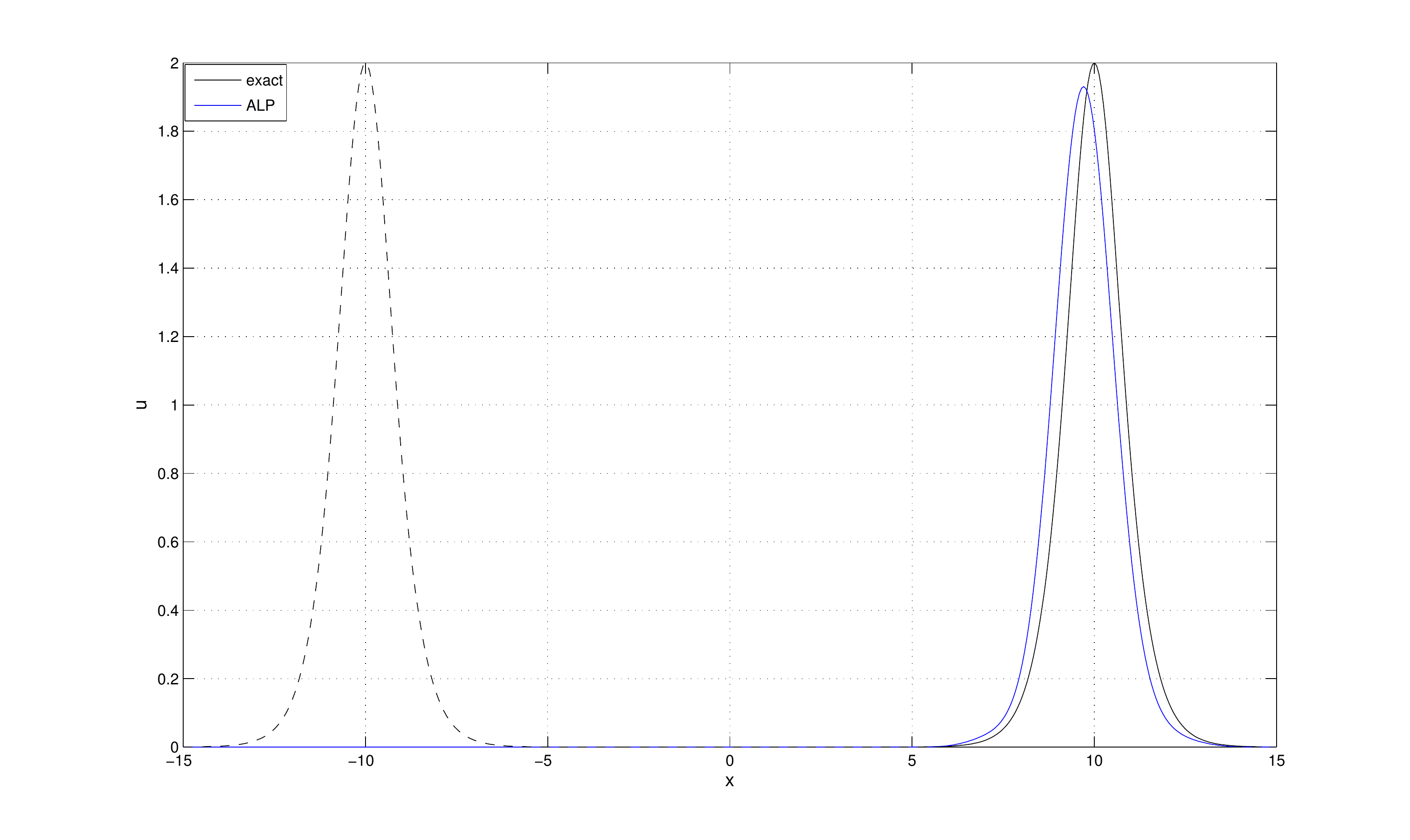}\\
\hspace{0.3cm} (a) & \hspace{0.3cm}(b) 
\end{tabular}}}
\caption{Comparison between the exact solution and the solution obtained with ALP when a) 10 modes and b) 25 modes are used in order to represent the evolution operator $\mathcal{M}(u)$. Dashed line represents the initial condition.}
\label{KdV_F2}
\end{figure}
In Fig.\ref{KdV_F2} a comparison between the exact (see Eq.(\ref{KdV_OneSolExact})) and the reduced-order integration method solutions is shown. In Fig.\ref{KdV_F2}.a) the final time solution is compared to the one obtained by representing the evolution operator $\mathcal{M}(u)$ with 10 modes:  the peak position, the amplitude and the wave shape are not correctly represented. Instead, when 25 modes are used, all the feature of the solution are well rendered, as it may be seen in Fig.\ref{KdV_F2}.b).

In Fig.\ref{KdV_F3}.b) the portrait of the operator $\mathcal{M}(u)$ at final time is shown. According to~\eqref{eq:mode_eq}, the entry $M_{mp}(u)$ represents the projection of the time derivative of the $m-$th mode onto the $p-$th one, at the current time. In this case, the highest values are all concentrated on the projection of the time evolution of the first mode (the only one used to represent $u$) on the others. Furthermore, the entries become smaller and smaller as $m,p$ increase, which is desirable when setting up a reduced order integration approach.

\begin{rem}
In this presentation, we restricted ourselves to P1 finite elements for the sake of conciseness. During our investigations, we also used P2 finite elements and Hermite polynomials, and we noticed that much less modes were necessary to achieve a good approximation in those cases. The link between the quality of representation by the modes and the quality of approximation of the underlying discretization method is an interesting question that might be addressed in future works. 
\end{rem}

\subsubsection{Three-soliton propagation}
\label{sec:three-solition}
\begin{figure}
\centerline{\hbox{\begin{tabular}{ccc}
\includegraphics[height=4.5cm,width=5.0cm]{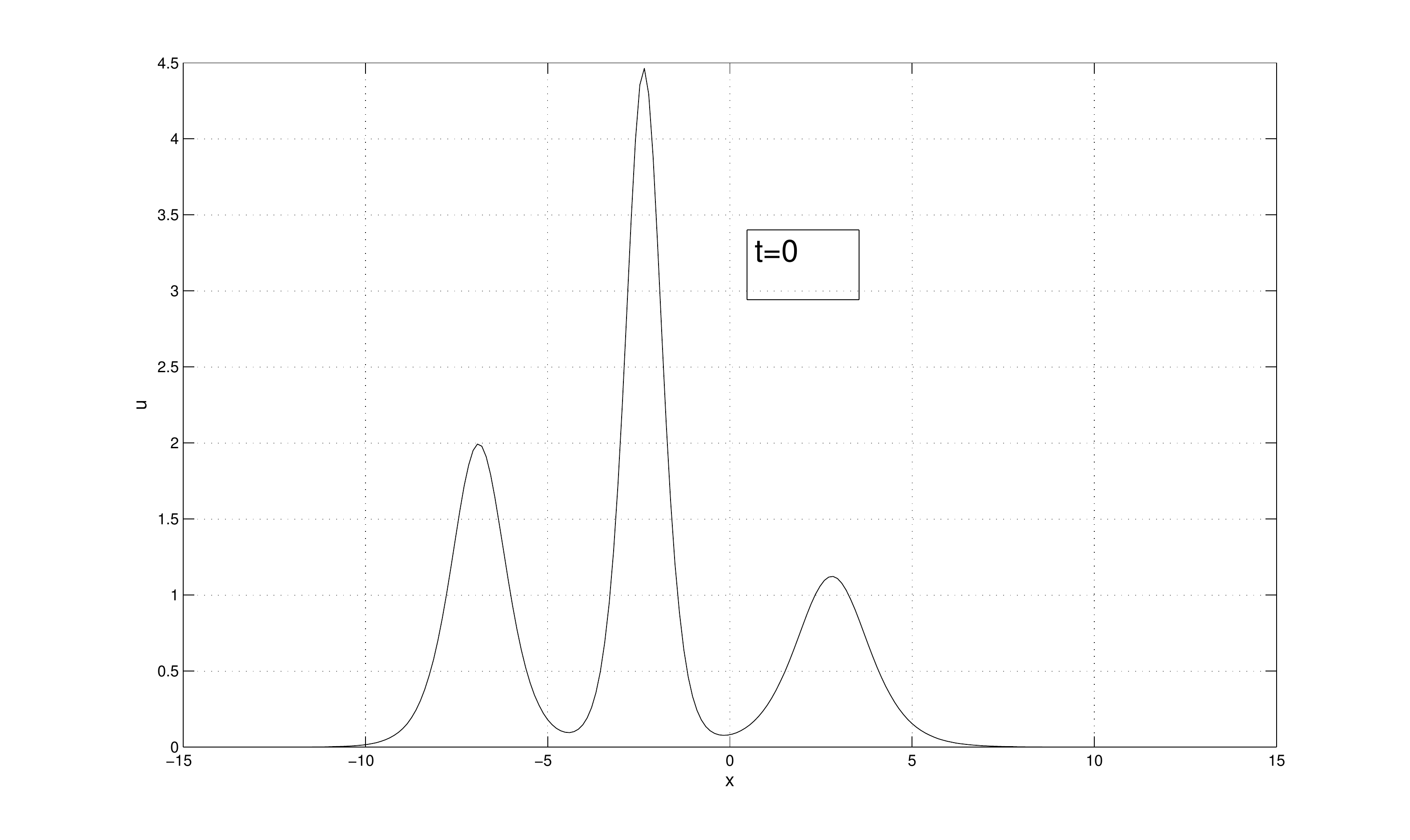} &
\includegraphics[height=4.5cm,width=5.0cm]{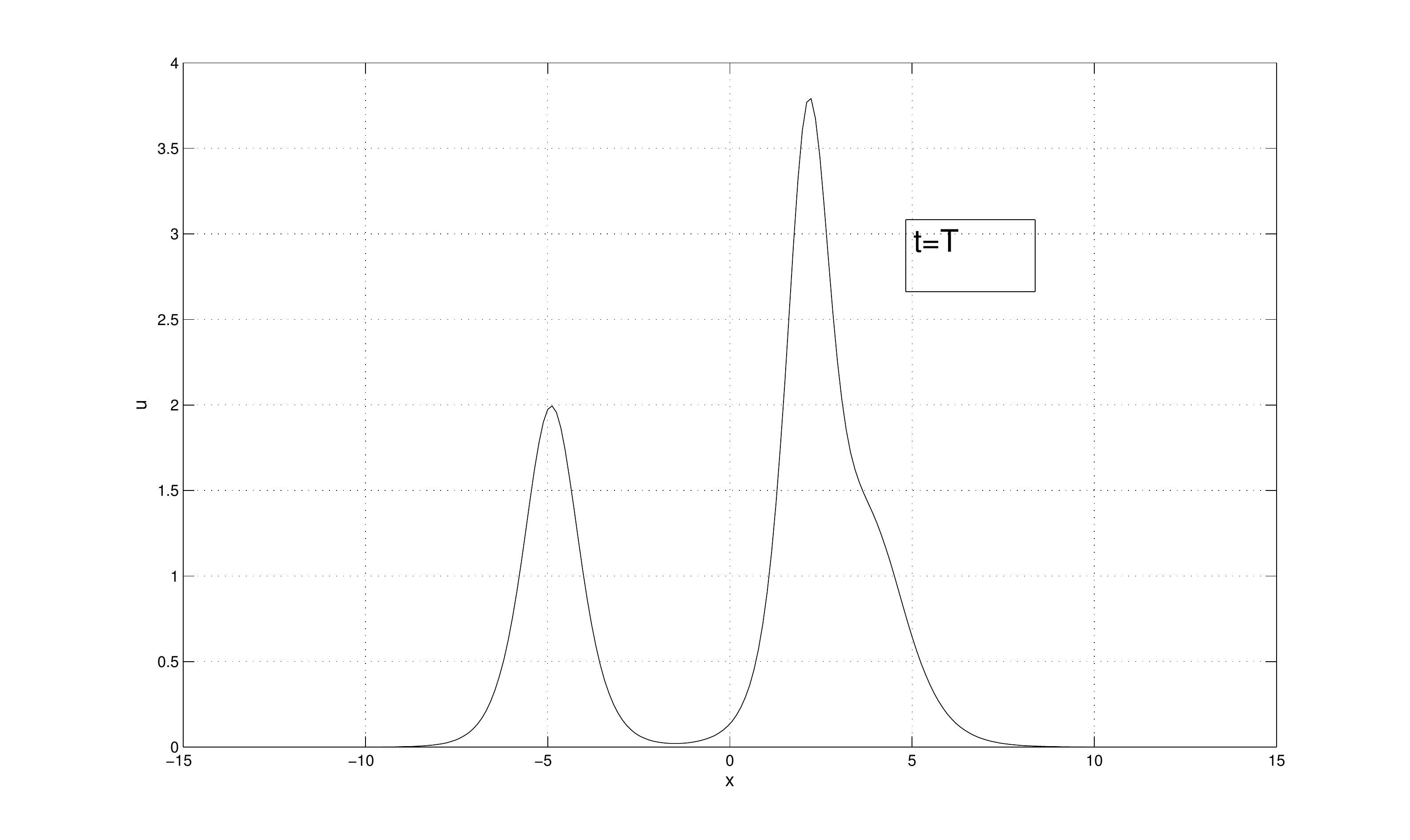} &
\includegraphics[height=4.5cm,width=5.0cm]{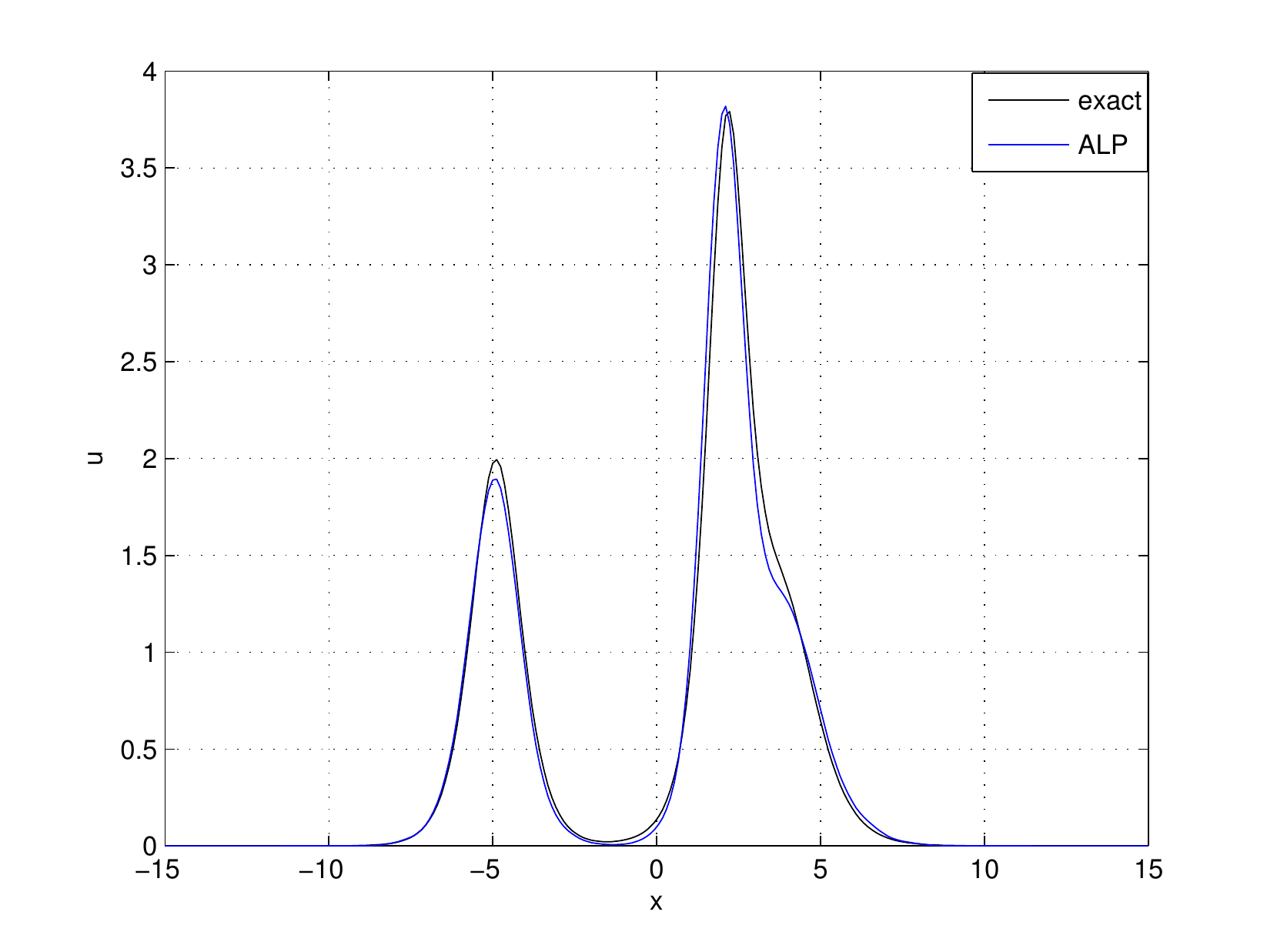}\\
\hspace{0.3cm} (a) & \hspace{0.3cm}(b) & \hspace{0.3cm}(c)
\end{tabular}}}
\caption{Plot of the exact three-soliton solution at a) initial and b) final times. In c), comparison between exact solution and reduced order model, at final time, when 20 modes are used to represent the operator. Problem setting defined by Eq.(\ref{threeSolExact}).}
\label{KdV_F4}
\end{figure}

In this section, we consider the propagation of a three-soliton by the KdV equation. The reference solution, shown in Fig.\ref{KdV_F4}.a-b), has been generated by considering the Gelfand-Levitan-Marchenko equation, that, when solved for the KdV equation (see \cite{Ablowitz_1}) provides:
\begin{equation}
u(x,t) = -2\partial^2_{x} \log (\det (I+A(x,t))),
\label{threeSolExact}
\end{equation} 
where $A\in\mathbb{R}^{n\times n}$ is the interacting matrix, written in terms of the scattering data. In particular:
\begin{equation}
A_{mn}(x,t) = \frac{c_m c_n}{k_m + k_n} \exp\left\lbrace (k_m+k_n)x - (k_m^3+k_n^3)t \right\rbrace,
\end{equation}
where $k_m,c_n$ are $2n$ scalar parameters that may be linked to position and speed of solitons (see \cite{Drazin_1,Shabat_1}).
For the present case: $c = [5.0\ 10^{-2},1.5\ 10^{-1},1.0\ 10^1]$, $k = [1.0, 1.5, 1.75]$, $x\in(-15,15)$ and $t\in(0,0.5)$. An interesting feature of this test case is the complex interaction between the waves: at the end of the simulation, two of them are fused together (see Fig.\ref{KdV_F4}.b)).

The scattering transform is applied at initial time and, by setting $\chi=1$, three distinct negative eigenvalues are found. This is in agreement with the analytical results and highlight the ability to decompose a traveling (nonlinearly interacting) waves system in modes which propagates separately. Thus, three modes are retained to represent the solution.

To represent the evolution operator $\mathcal{M}(u)$, we take $N_M=20$ modes. With this value of $N_M$, we checked at every time step that the relative difference in Frobenius norm $\mathcal{E}_{F}(N_M-1,N_M)$ (see \eqref{eq:error-frob}) was less than $1.0\%$. For the present case, the entries in the first three rows of its representation, corresponding to the projection of the time variation of the three waves used to recover $u$ on all the others, are higher in absolute value. It has been observed in the numerical experiments that, when $k$ modes are necessary and sufficient to have a good representation of $u(x,t)$, the matrix $M$ has the first $k$ rows (and columns, $M$ being skew-symmetric) which contains the entries that contribute, maximally, to the Frobenius norm of the operator. 

In Fig.\ref{KdV_F4}.b) a comparison at final time between the exact (see Eq.(\ref{threeSolExact})) and the ROM solutions is shown. The dynamics, in particular the coalescence phenomenon, is correctly represented.

\begin{figure}
\centerline{\hbox{\begin{tabular}{cc}
\includegraphics[height=5.0cm,width=6.5cm]{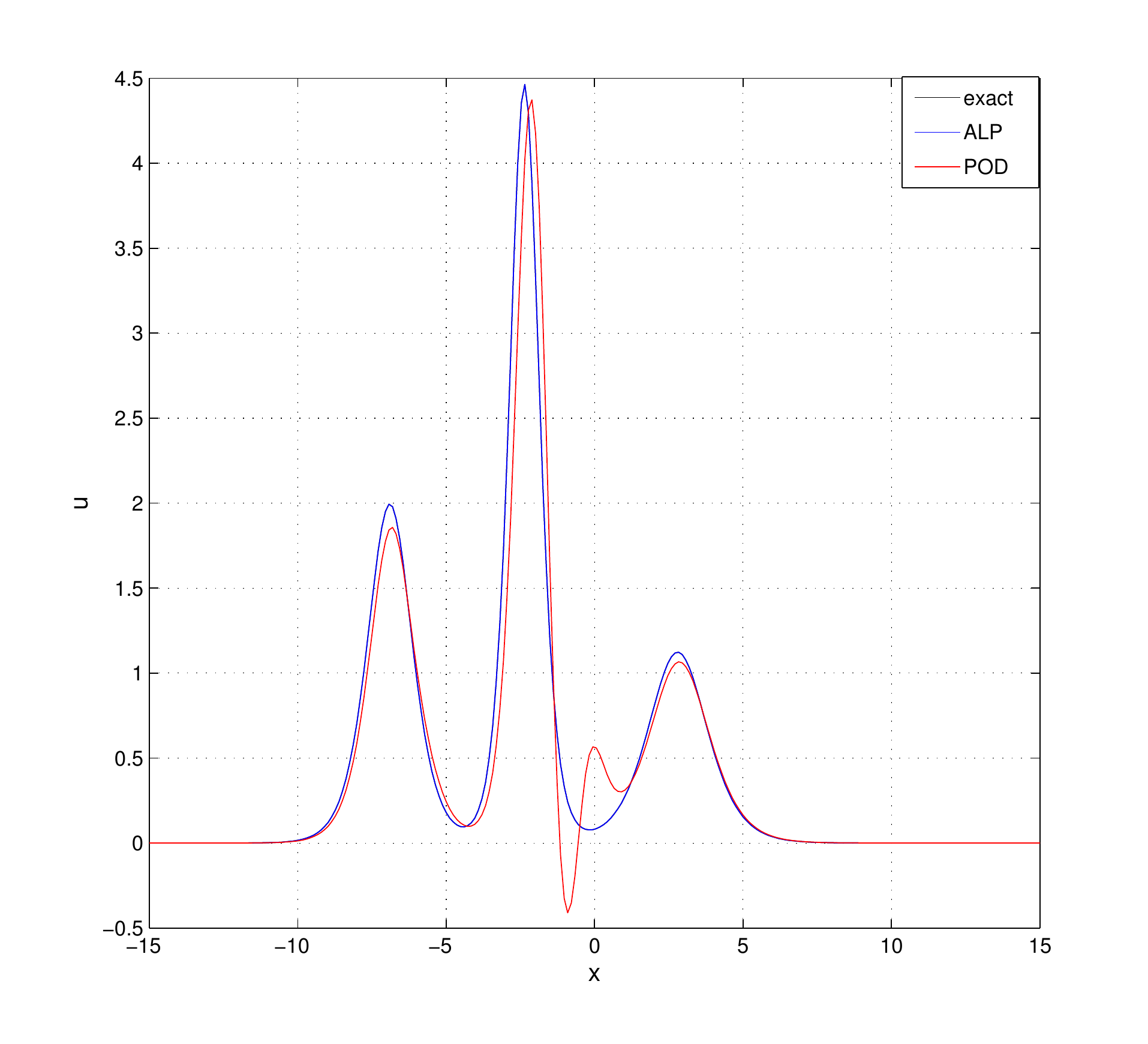} &
\includegraphics[height=5.0cm,width=6.5cm]{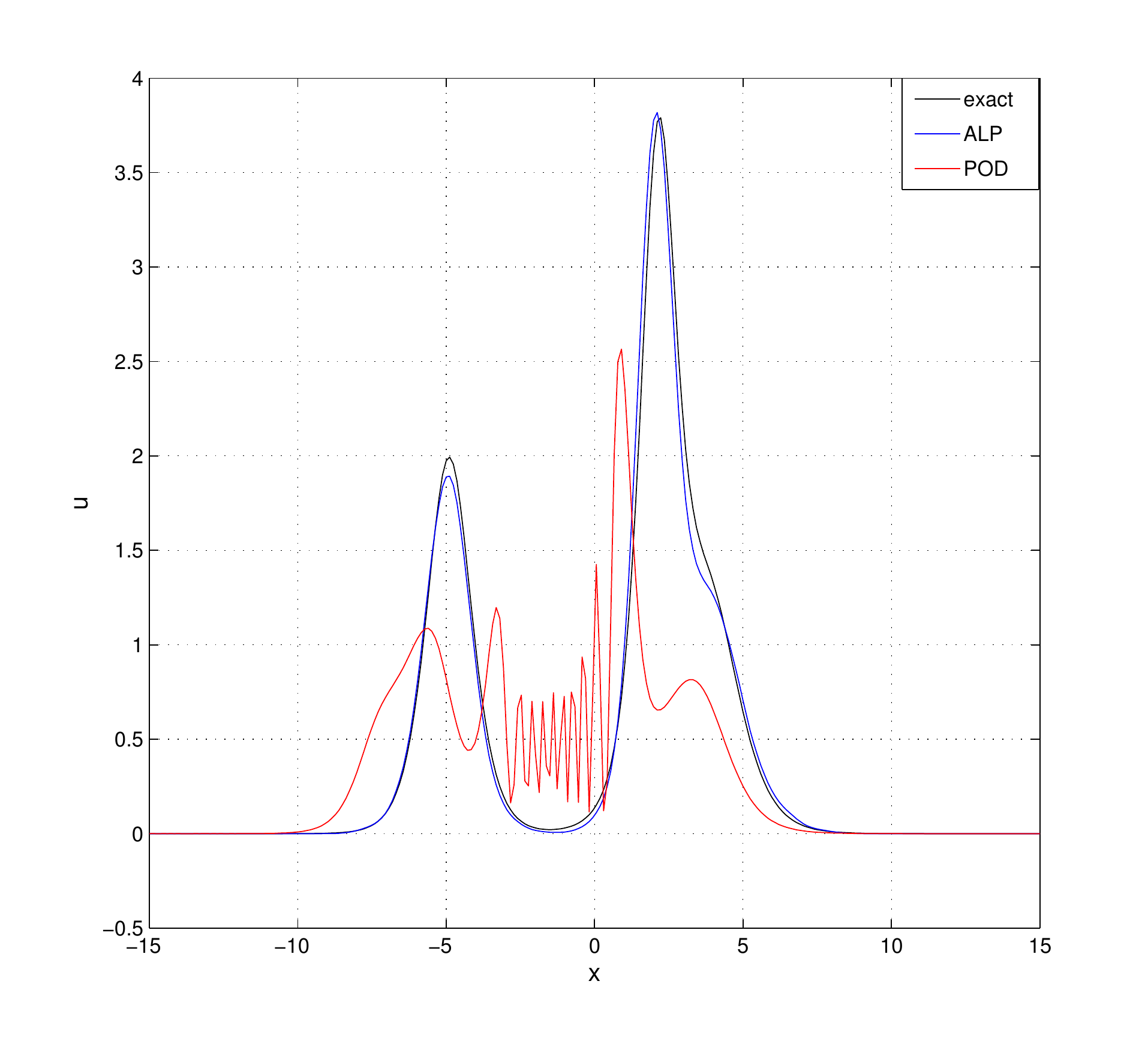}\\
\hspace{0.3cm} (a) & \hspace{0.3cm}(b) 
\end{tabular}}}
\caption{Comparison between POD and the ALP a) initial time, static reconstruction with 3 modes for both the techniques b) final time, integrated result with ALP, best possible reconstruction using 20 POD modes.}
\label{KdV_POD}
\end{figure}
Compared to the POD method (Proper Orthogonal Decomposition, see e.g. \cite{Sirovich_1}), an interesting feature of the ALP is its ability to extrapolate traveling wave type solutions out of the database used to defined the approximation space. To illustrate this point, we compared the two methods on the three-soliton test case. A POD basis was built by using the Sirovich technique (see \cite{Sirovich_1} for details), by taking $N_s = 50$ snapshots on half the time history, \emph{i.e.} between $t=0$ and $t=0.25$. In Fig.\ref{KdV_POD}.a), we compare the ability of the two approaches to approximate the initial data with only 3 modes. In Fig.\ref{KdV_POD}.b) we compare the solution at final time $t=0.5$. The POD reduced order model was run with 20 modes, and the ALP with $N_-=3$ modes to represent the solution and $N_M=20$ for $\mathcal{M}(u)$. In this example, the accuracy of the ALP clearly outperforms the POD. Let us emphasize that POD was built by using half of the time history while the ALP simply propagates the initial datum, without any pre-computed database.

\begin{figure}
\centerline{\hbox{\begin{tabular}{ccc}
\includegraphics[height=4.5cm,width=5.0cm]{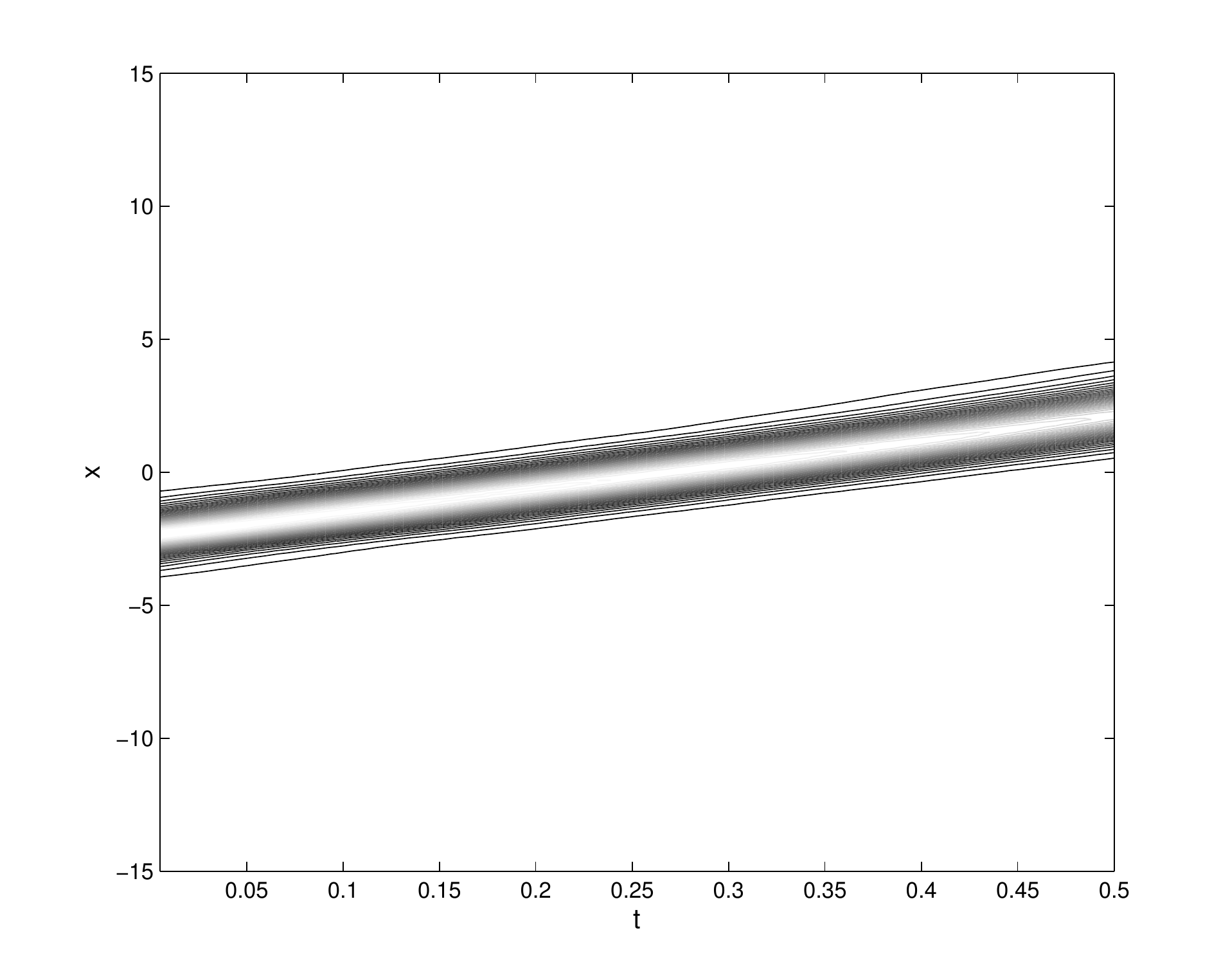} &
\includegraphics[height=4.5cm,width=5.0cm]{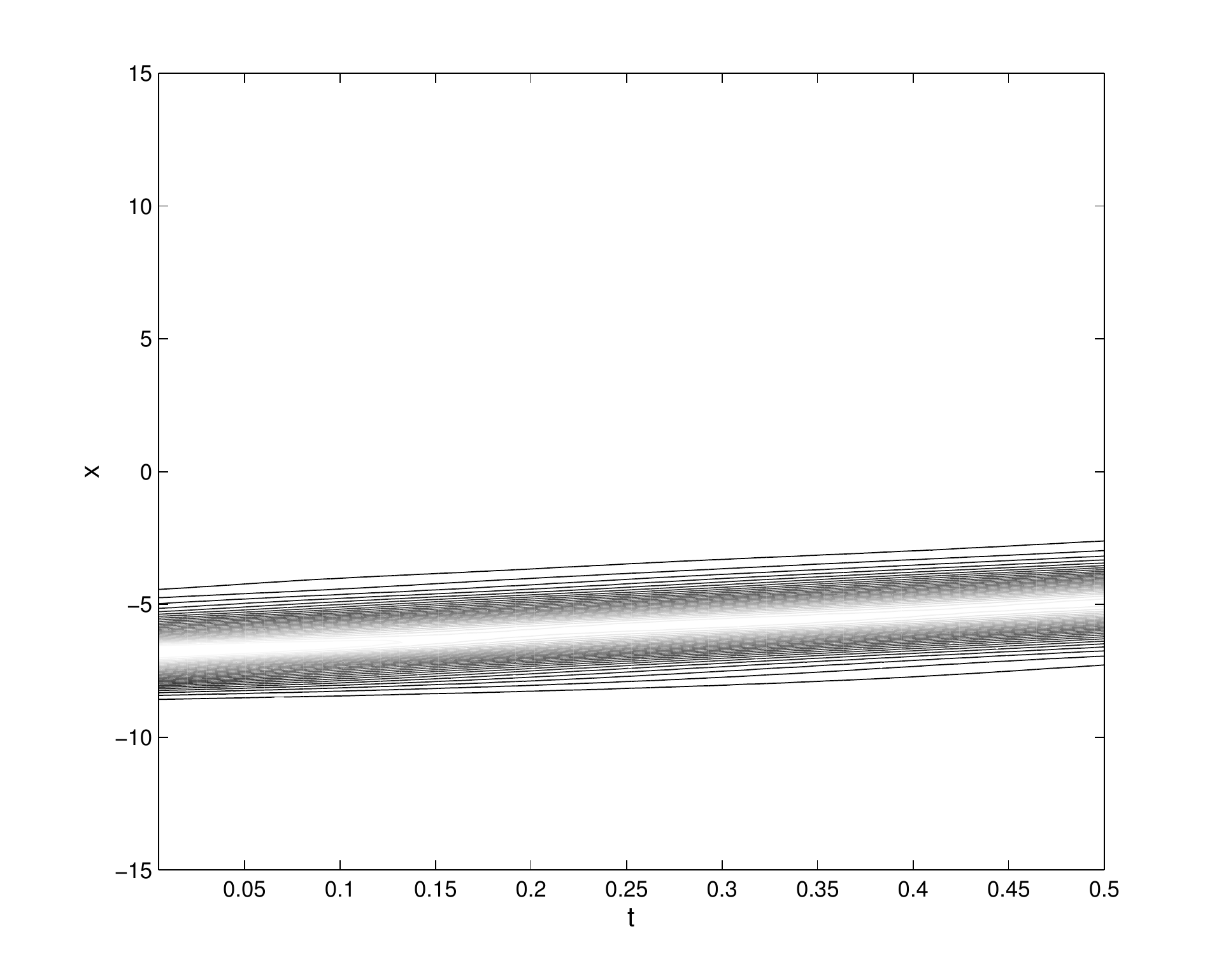}&
\includegraphics[height=4.5cm,width=5.0cm]{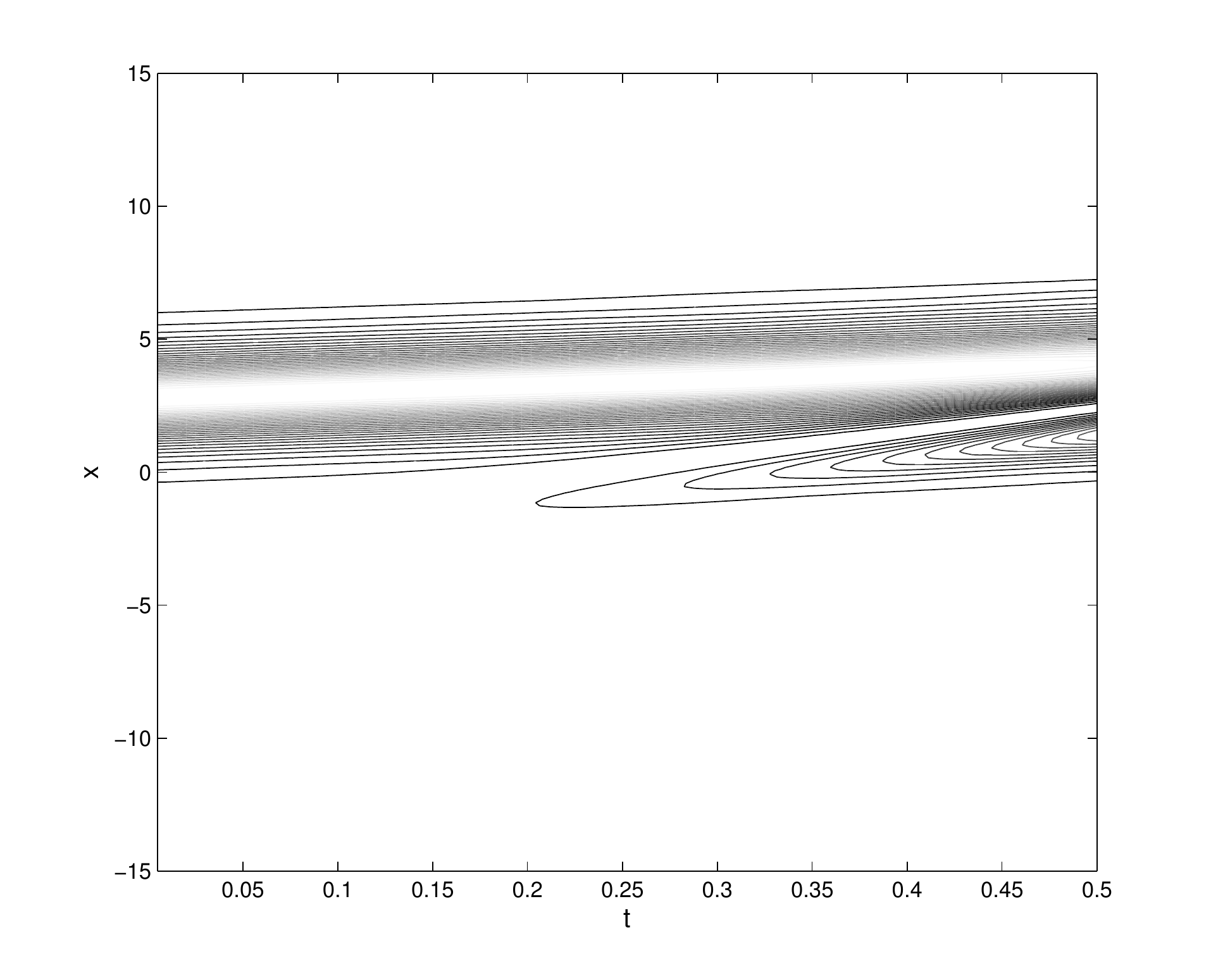}\\
\hspace{0.3cm} (a) & \hspace{0.3cm}(b) & \hspace{0.3cm}(c) 
\end{tabular}}}
\caption{Contour of the modes squared used to represent $u$ for a three-soliton propagation (see Eq.(\ref{threeSolExact})): a) first b) second c) third mode, 30 levels between maximum and minimum.}
\label{KdV_F5}
\end{figure}

In Fig.\ref{KdV_F5} the contour of the modes squared used to represent the solution -- i.e. the $\phi_m^2$ used in \eqref{eq:u-deift-trubo} -- is shown in $t,x$ plane. There are three distinct waves, propagated separately.
\begin{figure}
\centerline{\hbox{\begin{tabular}{ccc}
\includegraphics[height=4.5cm,width=5.0cm]{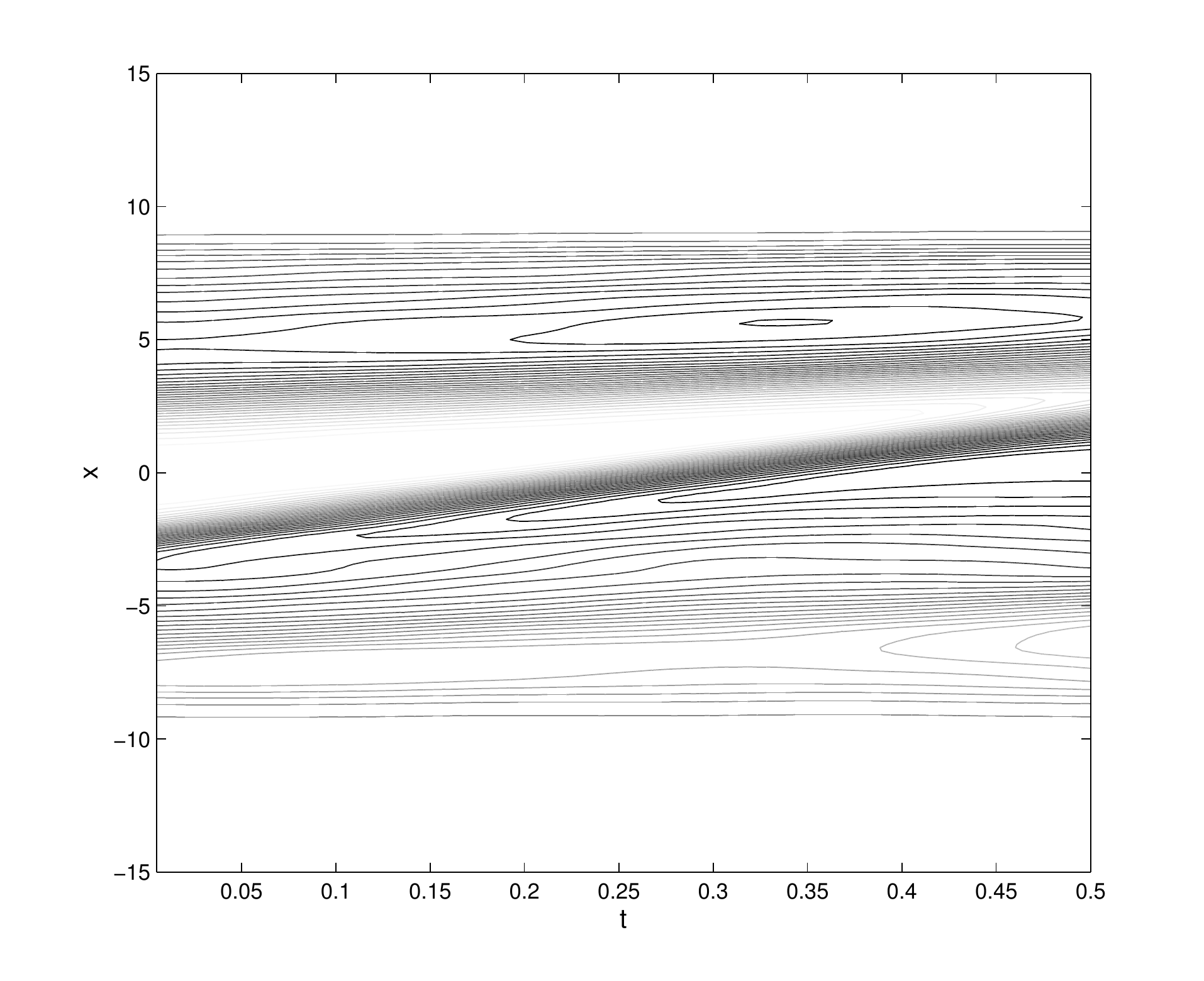} &
\includegraphics[height=4.5cm,width=5.0cm]{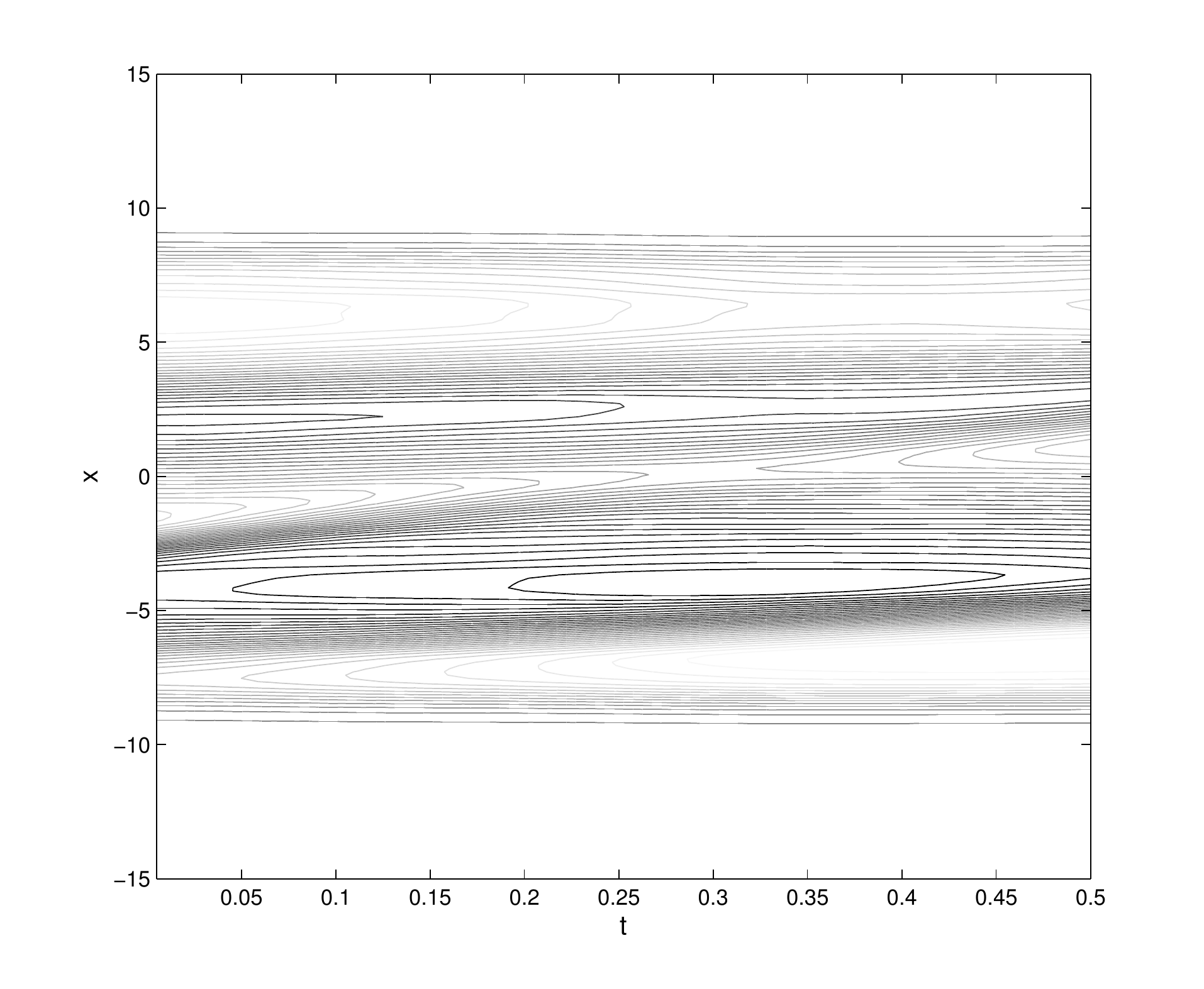}&
\includegraphics[height=4.5cm,width=5.0cm]{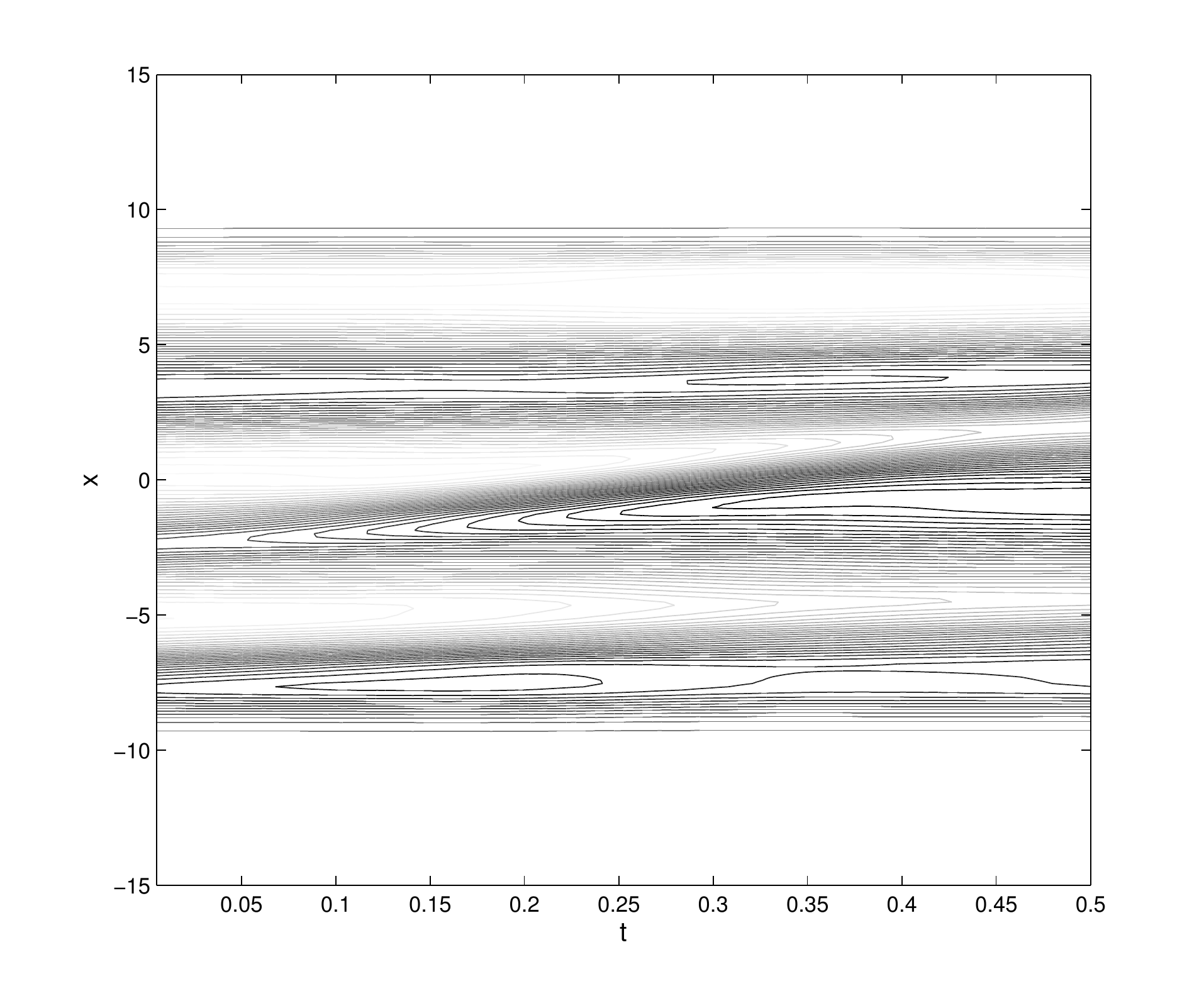}\\
\hspace{0.3cm} (a) & \hspace{0.3cm}(b) & \hspace{0.3cm}(c) 
\end{tabular}}}
\caption{Contour of the modes for a three-soliton propagation (see Eq.(\ref{threeSolExact})): a) fourth b) fith c) sixth mode, 30 levels between maximum and minimum.}
\label{KdV_F6}
\end{figure}
In Fig.\ref{KdV_F6} the contour of the first three modes $\psi_k$ associated with positive eigenvalues is shown. These modes are not used to reconstruct the solution, but they contribute to the waves evolution, in the sense that the projection of the time derivative of the waves on this modes is significant. 

\begin{rem}
	For the KdV equation, as well as for any equation admitting a Lax pair, the expression of the operator $\mathcal{M}(u)$ is known analytically. This may be profitably used for the reduced order model, since in that case,  there is no need to search for an approximated representation of the propagation operator. In other words, Steps 1 and 3 of the ALP algorithm can be replaced by the direct solution of equation~(\ref{ScattProb_1})$_2$. We checked for the KdV equation that the two approaches give similar results. 
\end{rem}

\subsection{Fisher-Kolmogorov-Petrovski-Piskunov equation}
In this section the FKPP equation is considered as an example of non-isospectral flow equation, in a finite domain, with homogeneous Dirichlet or Neumann boundary conditions. This differs substantially from the setting usually adopted when the Lax pair is known in a closed-form.



\subsubsection{1D FKPP with homogeneous Dirichlet boundary conditions}
\label{sec:fkpp-1d}

The equation reads:
\begin{equation}
\begin{split}
& \partial_t u = \partial_{xx}^2 u + \alpha u(1-u), \ in \ \Omega = [0,1]\\
& u = 0, \ on \ \partial \Omega
\end{split}
\label{FKPP_Eq1D}
\end{equation}
For the present case $\alpha = 1.0\cdot 10^3$, the space domain was $[0,1]$ and $N_f = 250$ piecewise linear functions were considered. The time domain was $[0,7.5\cdot 10^{-3}]$ and $100$ integration points were taken. The reference solution was obtained by discretizing in space by means of piecewise linear functions and by using a mixed implicit-explicit scheme in time: the linear diffusion part of the equation was discretized by means of a Cranck-Nicolson scheme, the nonlinear term by an explicit second order Adams-Bashforth scheme, with $\delta t = 7.5 \ 10^{-5}$.
\begin{figure}
\centerline{\hbox{\begin{tabular}{cc}
\includegraphics[height=5.0cm,width=6.5cm]{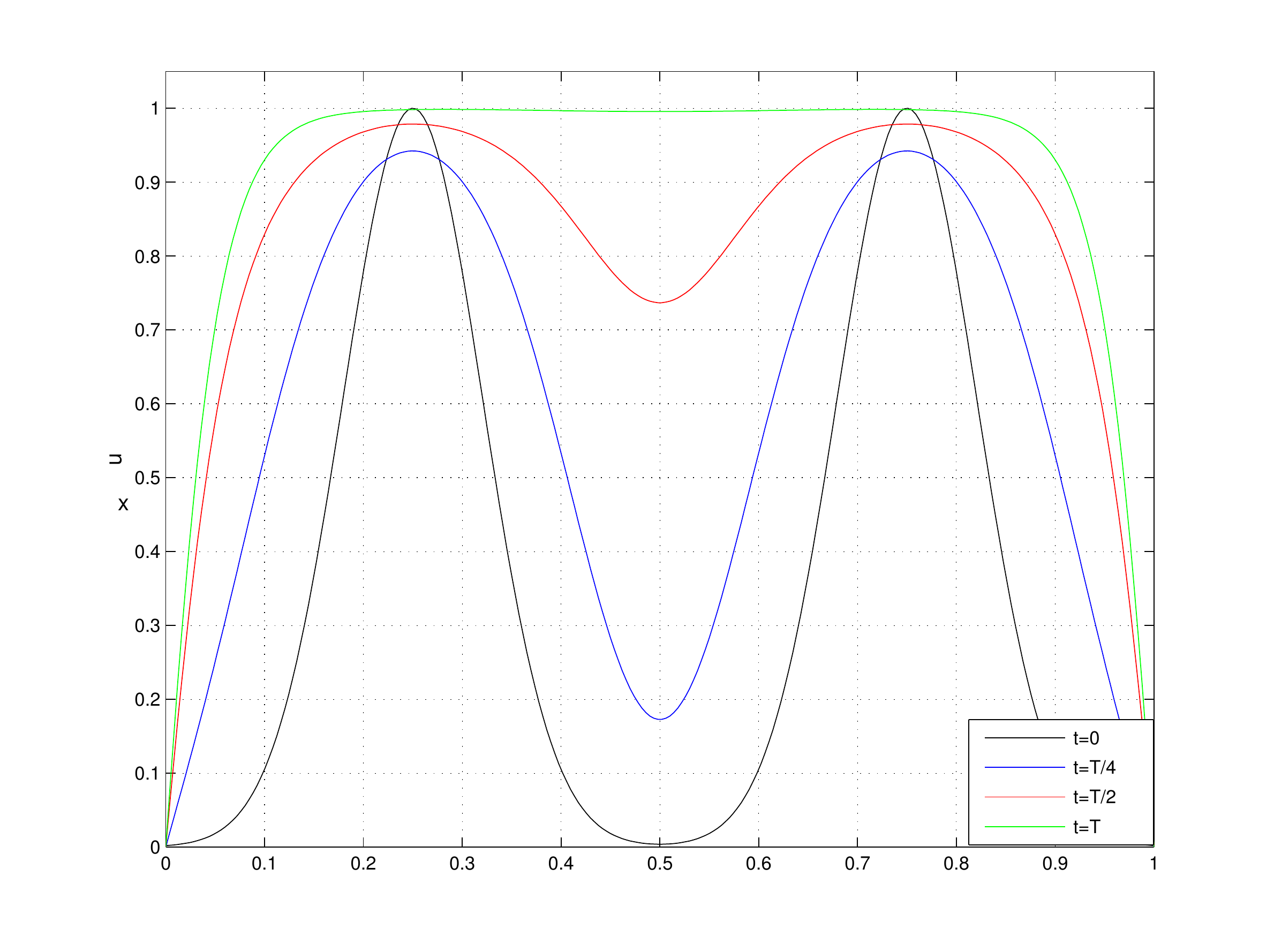} &
\includegraphics[height=5.0cm,width=6.5cm]{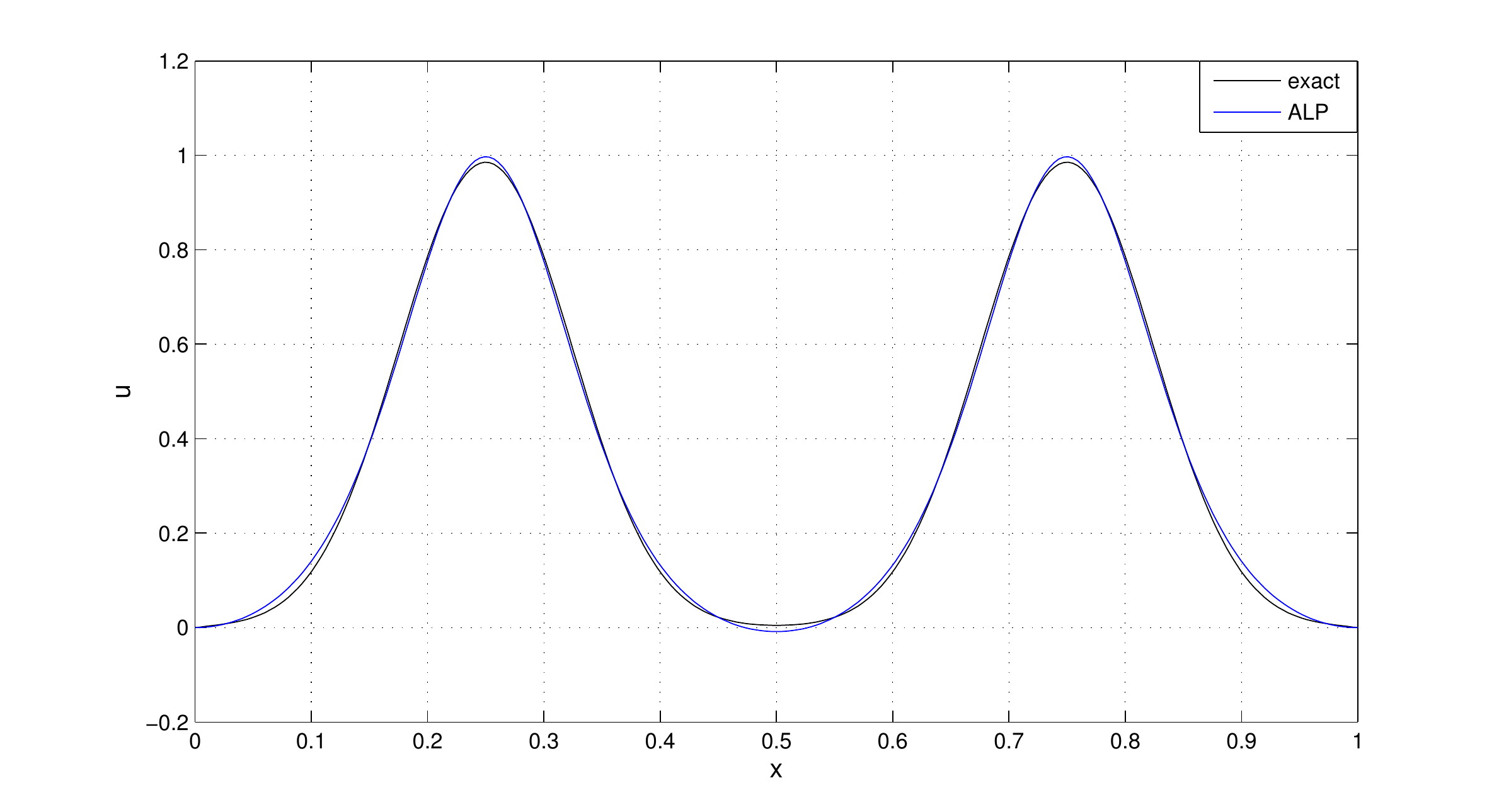}\\
\hspace{0.3cm} (a) & \hspace{0.3cm}(b) 
\end{tabular}}}
\caption{a) Plot of the reference solution of Eq.(\ref{FKPP_Eq1D}), at different times; b) Comparison between the exact solution and the reconstructed one at initial time by using four modes.}
\label{FKPP_Sol}
\end{figure}
In Fig.\ref{FKPP_Sol}.a) the reference solution is plotted at different times: the initial solution is taken as: $u_0 = \exp\left(-10^2(x-0.25)^2\right)+\exp\left(-100(x-0.75)^2\right)$. 

The number of modes to be retained in order to represent $u$ and $\mathcal{M}(u)$ have to be chosen. The first one is determined by setting $\chi$, for which no theoretical value is available, contrary to what happened for the KdV equation. The scattering parameter $\chi = 5\cdot10^2$ is chosen as explained in the initialization stage of the ALP algorithm, with $\epsilon_0=10^{-3}$. This corresponds to 4 distinct modes $\phi_k$ associated to negative eigenvalues. The number of $\psi_k$ modes for the approximation of $\mathcal{M}(u)$ is $N_M=10$, chosen with the same criterion used in the Section~\ref{sec:three-solition} (Frobenius norm modify by less that $1\%$ when adding a new mode).

In Fig.\ref{FKPP_Sol}.b) a comparison between the reference initial solution and its reconstruction is shown.
\begin{figure}
\centerline{\hbox{\begin{tabular}{cc}
\includegraphics[height=5.0cm,width=6.5cm]{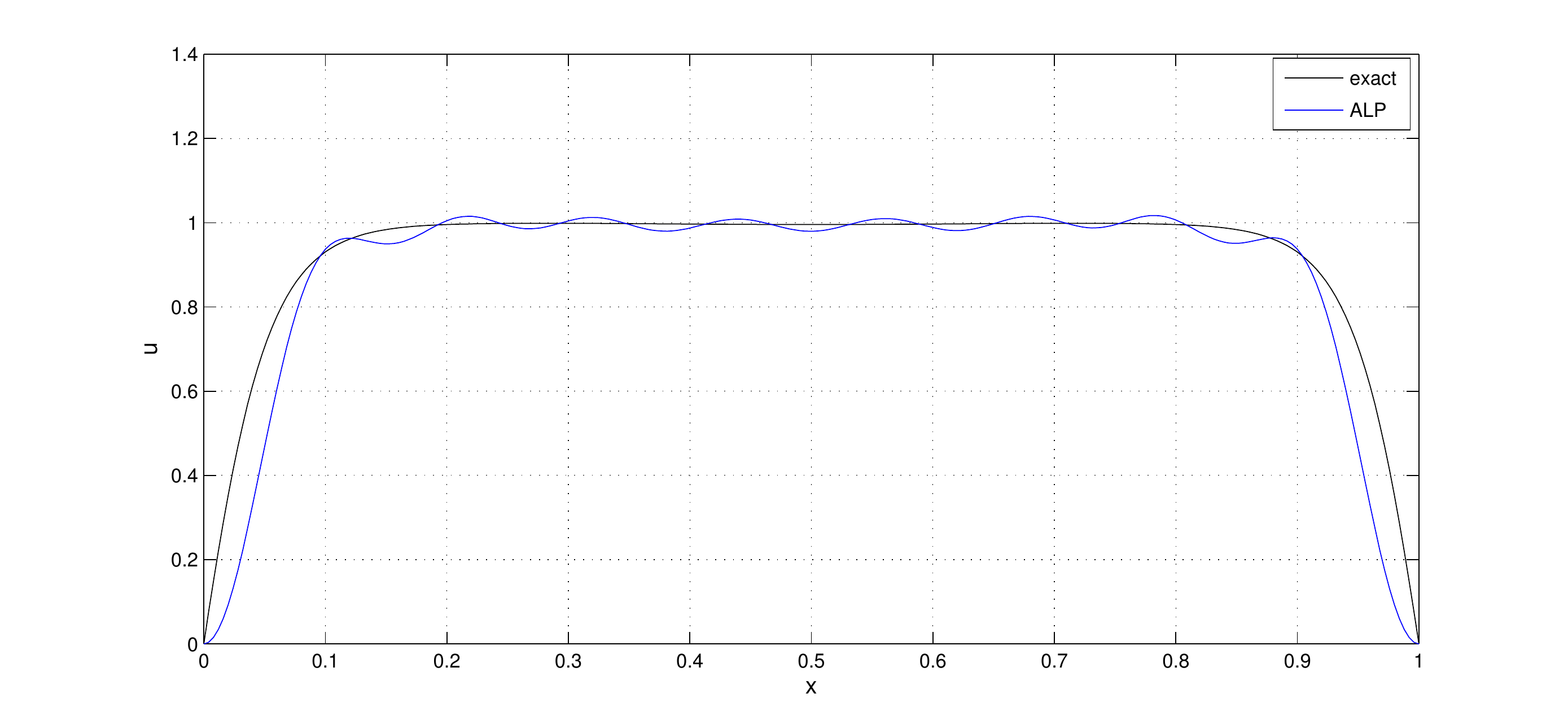} &
\includegraphics[height=5.0cm,width=6.5cm]{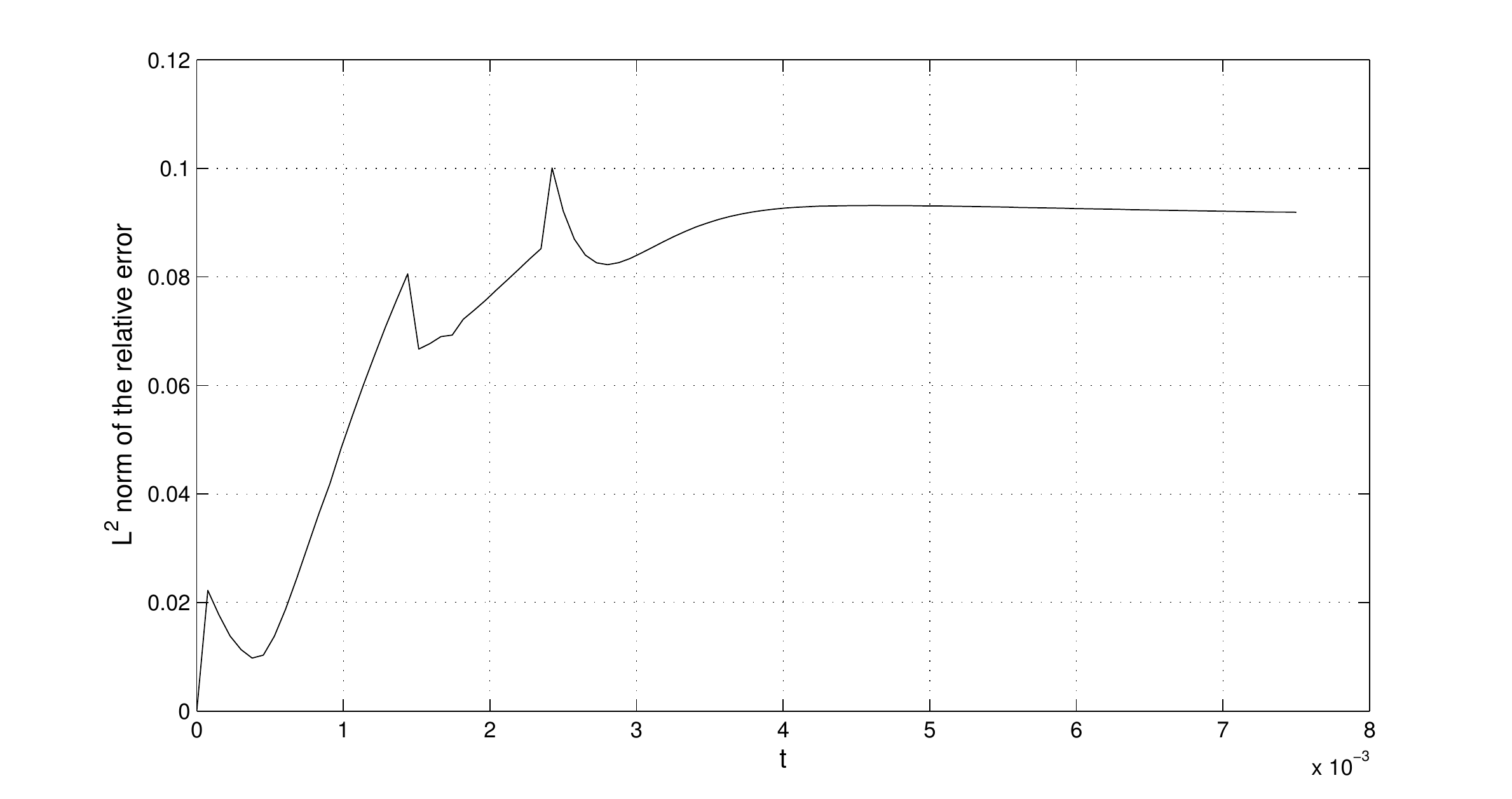}\\
\hspace{0.3cm} (a) & \hspace{0.3cm}(b) 
\end{tabular}}}
\caption{a) Comparison between the exact solution of Eq.(\ref{FKPP_Eq1D}) and the reconstructed one at final time; b) relative error in $L^2$ norm with respect to time}
\label{FKPP_Err}
\end{figure} 
 
In Fig.\ref{FKPP_Err}.a) a comparison between the reference solution at $t=T$ and the reconstructed one is shown. The error is mainly due to the appearence of wavy structures. However, the dynamics, wich is non-trivial (coalescence of structures) is correctly recovered. 

In Fig.\ref{FKPP_Err}.b) the relative error in $L^2$ norm is shown as a function of time. It is smaller than $10\%$, and in this case it is not monotonically increasing in time. This behavior is due to the eigenvalues dynamics. Contrary to what happened with the KdV equations, the eigenvalues are evolving in time and some positive eigenvalues may become negative.  At the very beginning, $N_-=4$ modes (corresponding to the negative eigenvalues) are sufficient to represent the solution. Although it is possible to keep these modes throughout the simulation, we noticed that the approximation was improved if, during the simulation, we added the modes corresponding to the new negative eigenvalues, as indicated in Remark~\ref{rem:variable-num-mode}.
In Fig.\ref{FKPP_Err}.b) the error peaks at times $t = 1.5\cdot 10^{-3}$ and $t = 2.5\cdot 10^{-3}$ corresponds to the instant at which $\lambda_5$ and $\lambda_6$ respectively change their sign, and the corresponding $\phi_5$ and $\phi_6$ are used to approximate the solution. 

\subsubsection{2D FKPP with homogeneous Neumann boundary conditions}

\begin{figure}
\centerline{\hbox{\begin{tabular}{ccc}
\includegraphics[height=3.5cm,width=4.0cm]{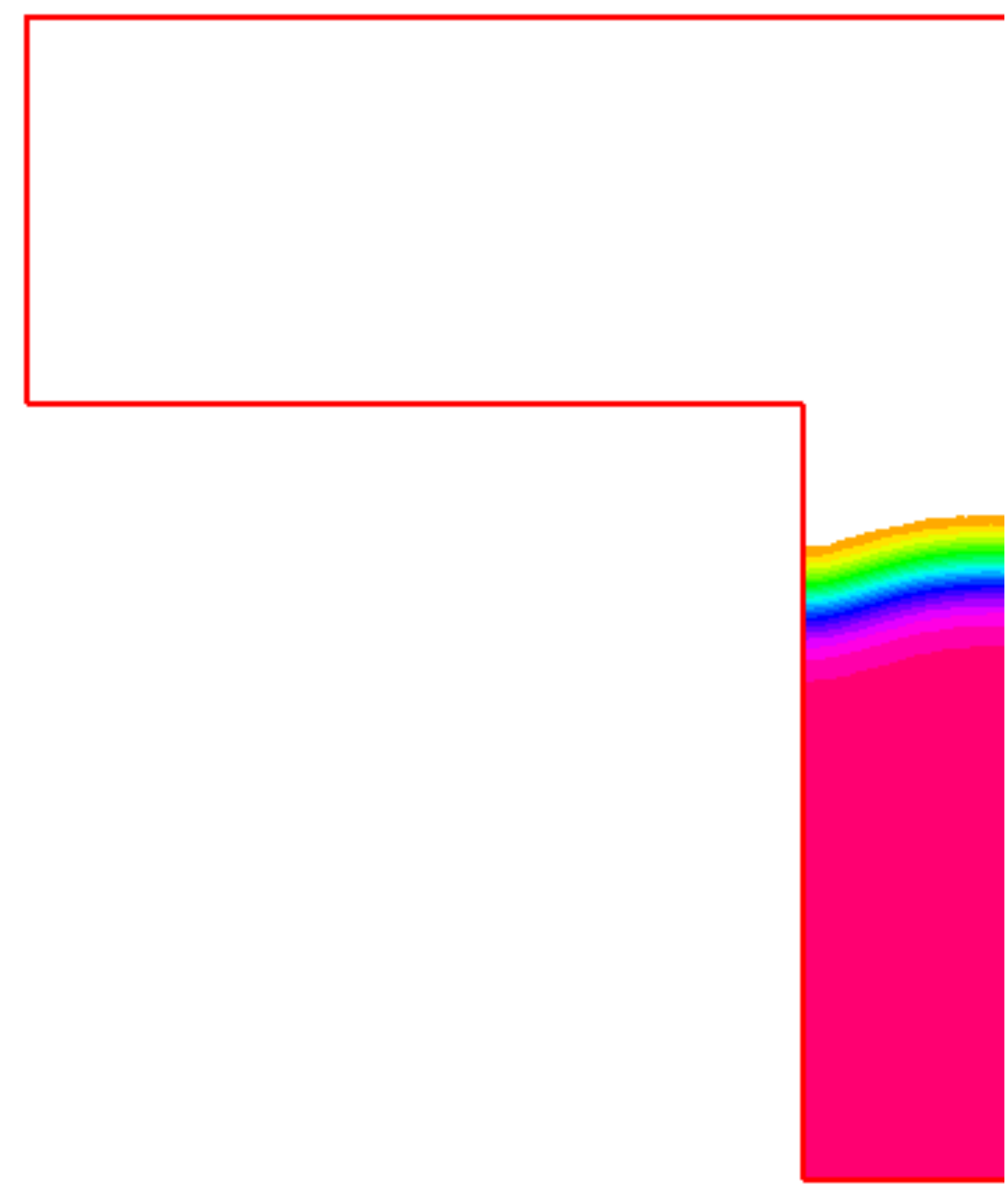} &
\includegraphics[height=3.5cm,width=4.0cm]{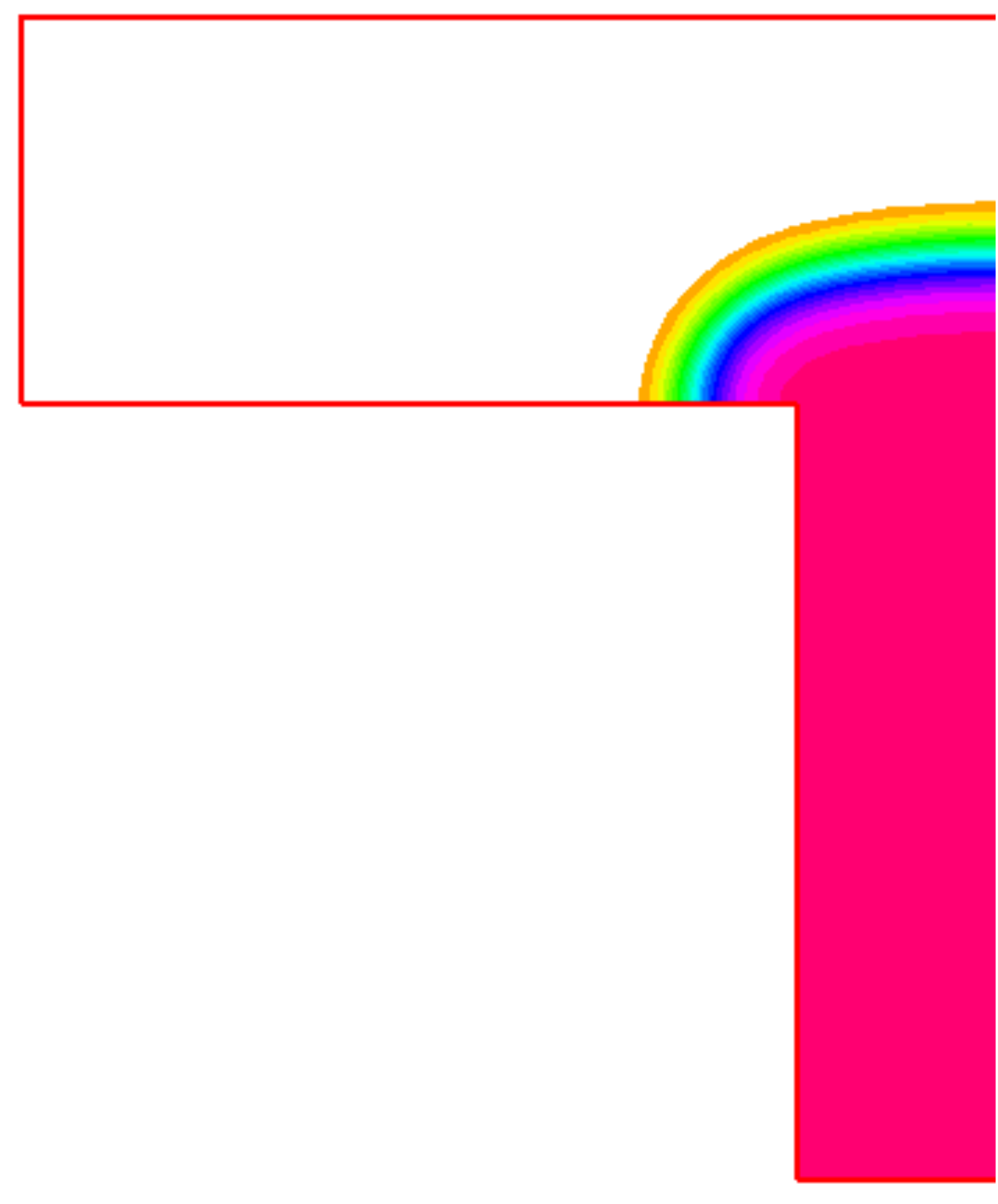} &
\includegraphics[height=3.5cm,width=4.0cm]{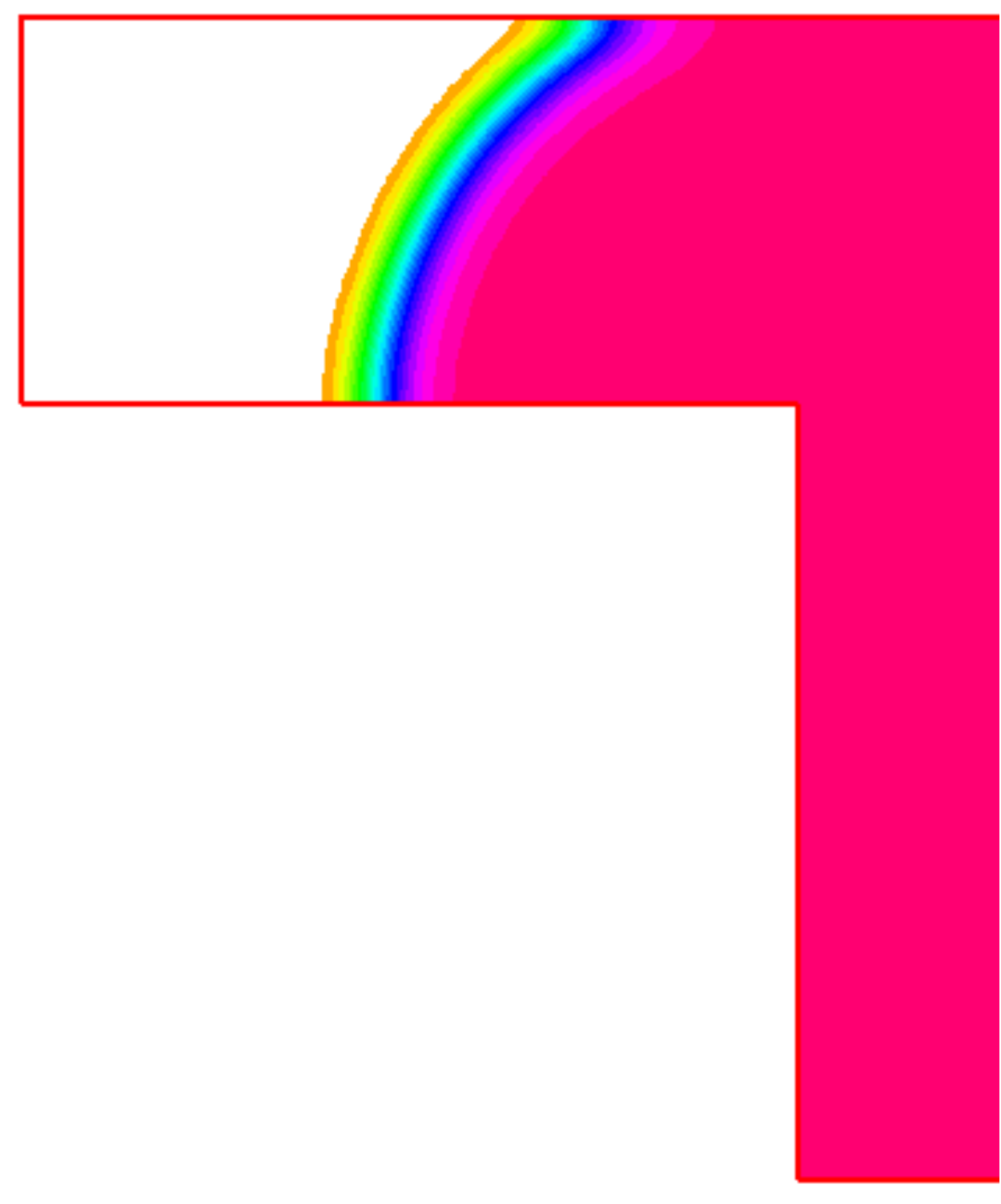}\\
\hspace{0.2cm} (a) & \hspace{0.2cm}(b) & \hspace{0.2cm} (c)
\end{tabular}}}
\caption{Contour levels of the reference solution of Eq.(\ref{FKPP_Eq2D}) on half the domain, 30 levels between maximum (red) and minimum (orange) at a) t=0, b) t=T/2, c=T.}
\label{FKPP2D_refSol}
\end{figure}

The 2D FKPP equation reads:
\begin{equation}
\begin{split}
\partial_t u = \Delta u + \alpha u (1-u), \ \ in \ \Omega,
\partial_n u = 0 \ \ on \ \partial\Omega,
\end{split}
\label{FKPP_Eq2D}
\end{equation}
where $\alpha = 10^3$ for the presented test case and homogeneous Neumann boundary conditions are imposed on $\partial\Omega, \ \Omega\subseteq \mathbb{R}^2$. 
\begin{figure}
\centerline{\hbox{\begin{tabular}{cc}
\includegraphics[height=5.0cm,width=6.5cm]{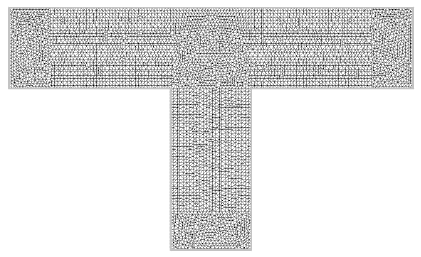} &
\includegraphics[height=5.0cm,width=6.5cm]{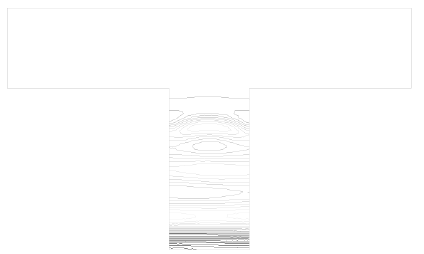}\\
\hspace{0.3cm} (a) & \hspace{0.3cm}(b) 
\end{tabular}}}
\caption{a) Geometry of the domain $\Omega$ and related mesh for Eq.(\ref{FKPP_Eq2D}) discretization b) Contour for the relative error field on the initial datum approximation, 25 isoline between maximum and minimum.}
\label{FKPP2D_Geom}
\end{figure}
A reference solution for this problem was obtained by integrating the equation by means of a P1 finite element method. The geometry of the domain and the mesh are shown in Fig.\ref{FKPP2D_Geom}.a). The boundary is a set of line segments of unitary and integer multiple of unity length; the mesh counts $9490$ triangles. In order to advance the initial solution in time a hybrid Crank-Nicolson, Adams-Bashforth 2 scheme was implemented, to discretize the diffusive and the reaction terms respectively. The time step was chosen as $\delta t = 2.5\ 10^{-4}$. The initial condition (see Fig.\ref{FKPP2D_refSol}.a) ) was obtained by making a Gaussian $u_g = \exp(-50((x-2.5)^2 + (y-0.5)^2)$ evolve for $n=75$ time steps. 
The reference solution is featured by a nonlinear front propagation. In particular, this geometry makes the evolution of this system particularly challenging to be recovered. First, the front propagates upward, then, once arrived at the $T$ bifurcation it starts spreading almost spherically (see Fig.\ref{FKPP2D_refSol}.b)), and finally, due to Neumann boundary condition, when it reaches the upper wall, it splits and two fronts start propagating horizontally with $\pm x$ direction (see Fig.\ref{FKPP2D_refSol}.c) ). The overall evolution took $n=125$ iterations. 

The initial datum was recovered by means of $N_{-}=6$ modes, obtained by solving the scattering problem with $\chi = 40$. In Fig.\ref{FKPP2D_Geom}.b) the contour of the initial relative error field is represented. It is featured by oscillations, especially located near the lower bound of the domain, and it is of overall $10\%$ in $L^2$ norm. It is due to the fact that the laplacian operator is not suitable for this kind of equation in finite domain, and the use of another compact operator should improve the quality of the datum reconstruction. The choice of the regularization operator is essentially a guess even in an analytical approach to the scattering transform (see \cite{Ablowitz_1} for a detailed discussion on this topic). This point will be the object of further works.

\begin{figure}
\centerline{\hbox{\begin{tabular}{ccc}
\includegraphics[height=3.5cm,width=4.0cm]{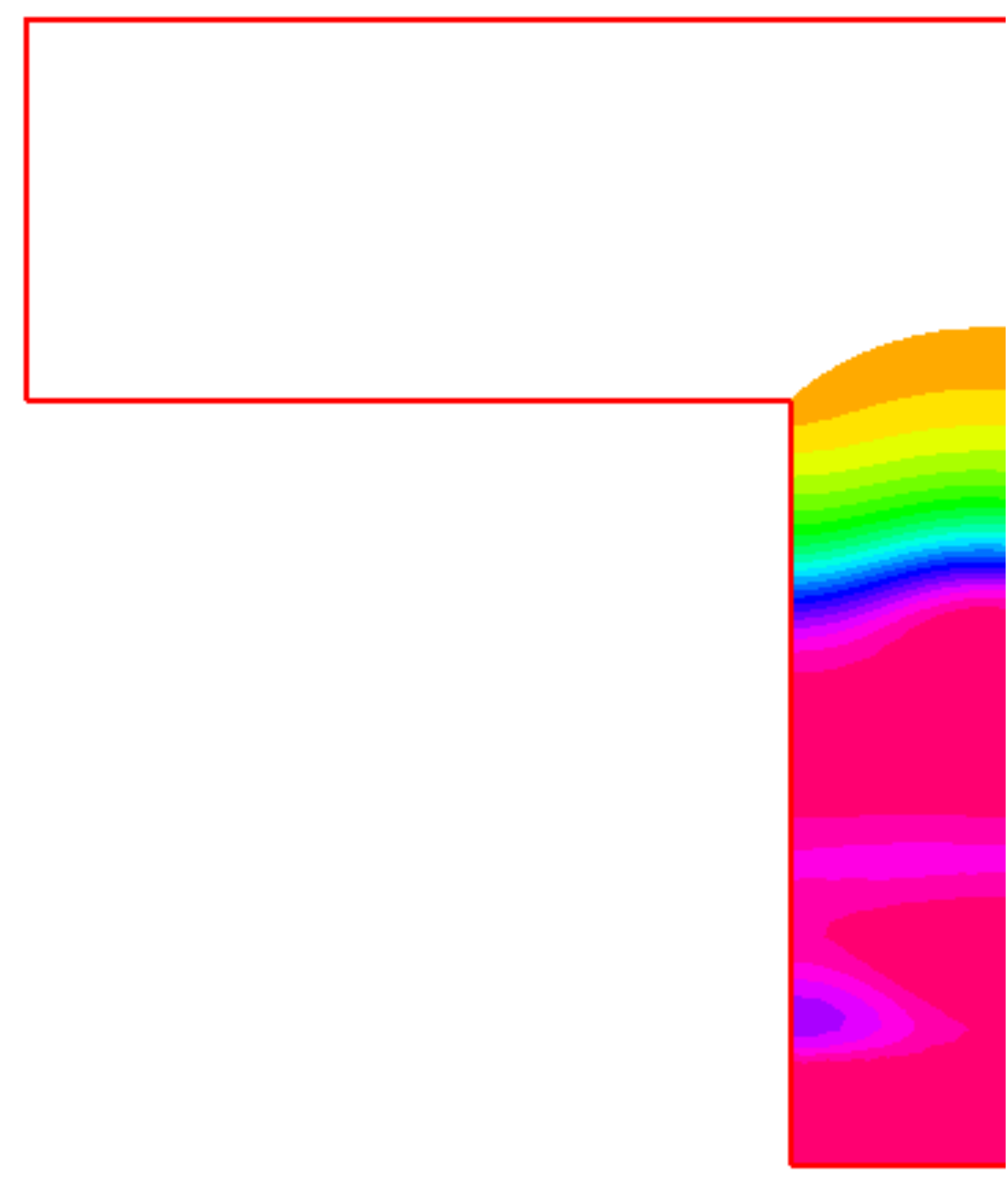} &
\includegraphics[height=3.5cm,width=4.0cm]{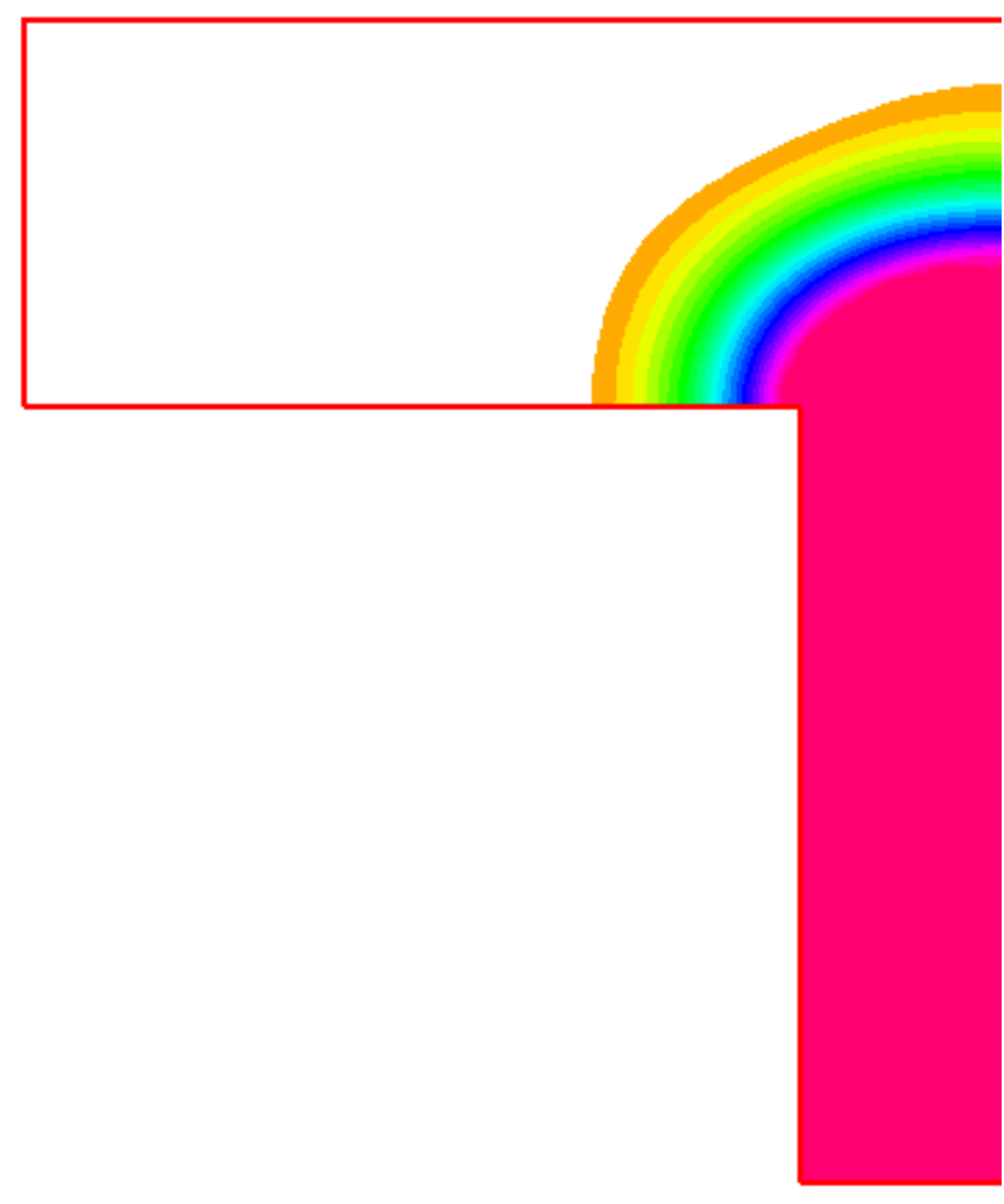} &
\includegraphics[height=3.5cm,width=4.0cm]{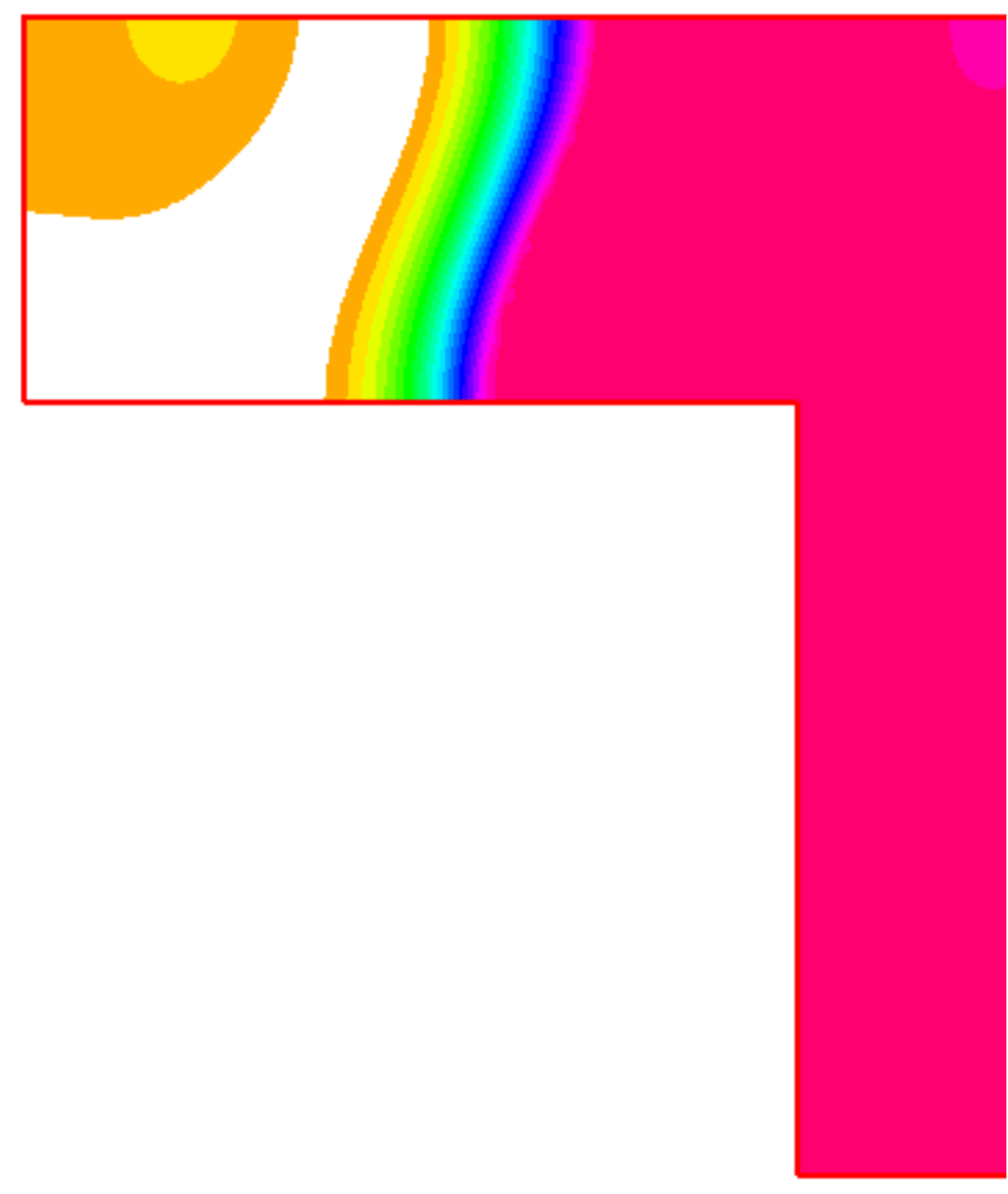}\\
\hspace{0.2cm} (a) & \hspace{0.2cm}(b) & \hspace{0.2cm} (c)
\end{tabular}}}
\caption{Contour levels of the solution obtained by ALP, visualized on half the domain, 30 levels between maximum (red) and minimum (orange) at a) t=0, b) t=T/2, c=T.}
\label{FKPP2D_simSol}
\end{figure}
In Fig.\ref{FKPP2D_simSol} the contour levels of the solution obtained by propagating the modes are shown, at the same time instants of those ones shown for the reference solution, using the same scale. Despite the fact that the initial datum reconstruction is affected by some error, all the elements of the dynamics, that is the front speed and the global movement are recovered in a satisfactory way. Some improvement is needed on the accuracy in the solution representation (\emph{i.e.} the front shape). 

This represents another case in wich POD would not be able to extrapolate out of database (by construction) and tends to perform poorly in terms of reduction.

The modes were  propagated by means of the proposed algorithm with a time step of $\delta t = 1.0 \ 10^{-3}$.
\begin{figure}
\centerline{\hbox{\begin{tabular}{ccc}
\includegraphics[height=3.5cm,width=4.0cm]{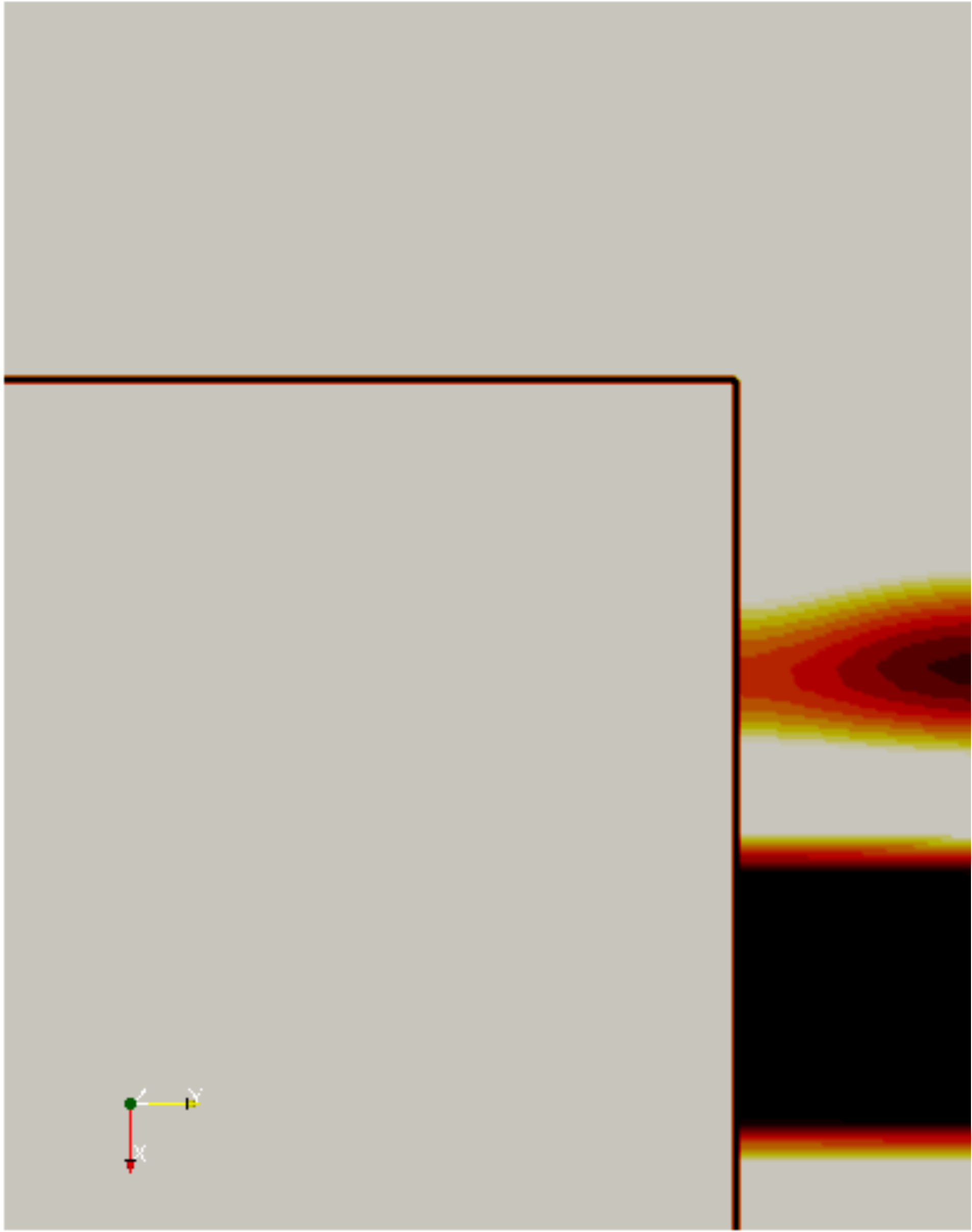} &
\includegraphics[height=3.5cm,width=4.0cm]{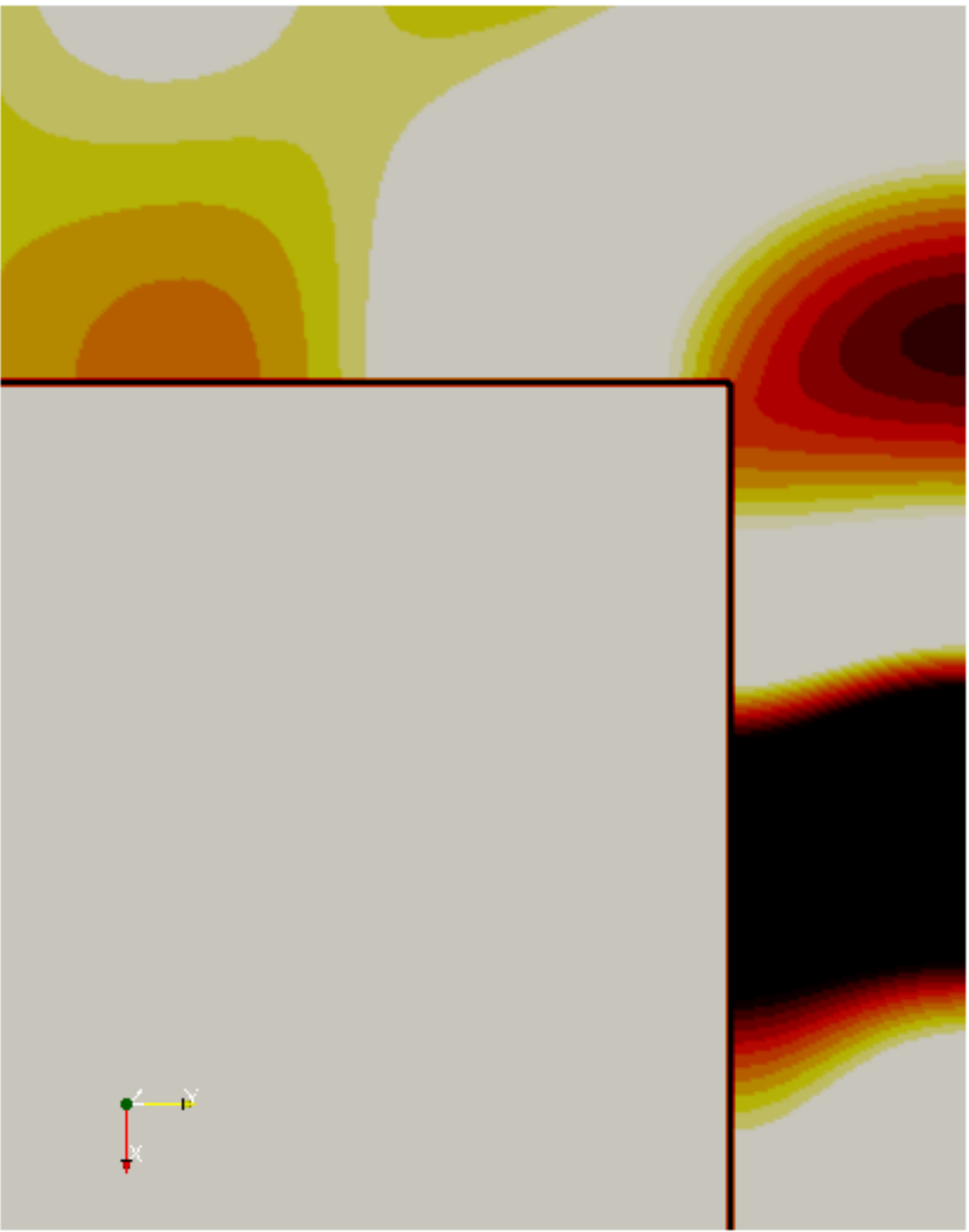} &
\includegraphics[height=3.5cm,width=4.0cm]{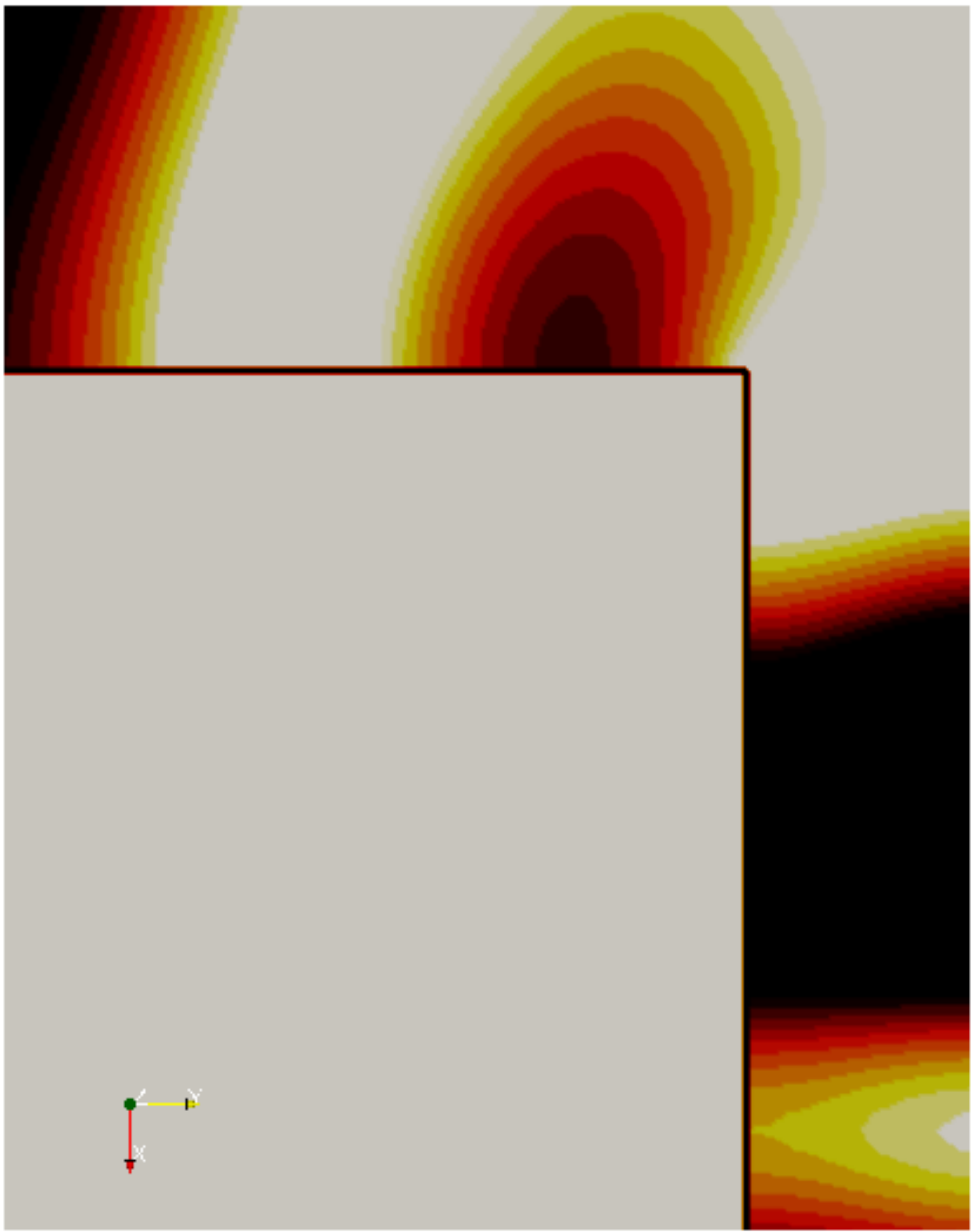}\\
\hspace{0.2cm} (a) & \hspace{0.2cm}(b) & \hspace{0.2cm} (c)
\end{tabular}}}
\caption{Contour levels of the most significant mode (fith one) on half the domain, 30 isolines between maximum and minimum at a) t=0, b) t=T/2, c=T.}
\label{FKPP2D_Mode4}
\end{figure}

In Fig.\ref{FKPP2D_Mode4}.a-c) the contour of the most significant mode is represented, at three different times, on half the domain (it is symmetric with respect to the $y$ axis). The mode travels and then splits, its dynamics being similar to the one of the solution, featured by large displacements.

\section{Conclusions and perspectives}

We have proposed a new reduced-order model technique, called ALP, consisting of three stages. First, a set of modes is constructed by the Semi Classical Signal Analysis method, i.e. as a set of orthonormal eigenfunctions of a linear Schr\"odinger operator associated to the initial condition (approximated scattering transform). Second, the modes are propagated by an approximation of a Lax operator, and the eigenvalues are updated. Third, the solution is reconstructed by solving a problem that plays a role similar to the inverse scattering transform. 

This approach allows us to set up a reduced order discretization that seems to be efficient for those problems involving progressive wave or front propagation. Contrary to other reduced-order methods, like POD, the solution can be extrapolated out of the database of pre-computed solutions. The method has been successfully tested on the KdV and FKPP equations, in 1D and 2D. 

The application of ALP to other problems is currently under investigation, in particular to a set of Euler equations modeling a network of arteries and to cardiac electrophysiology problems. 

Many questions would deserve further investigations: the number of modes used to approximate the solution or the propagation operator could be adapted along the resolution, for example based on the indicator \eqref{eq:error-frob}; other operators than the Laplacian might used for the scattering step; other schemes could be used to solve \eqref{eq:psi-tilde}; positive eigenvalues might be used to approximate the solution, etc. This will be the subject of future works.

\bibliography{ROMBiblio}
\bibliographystyle{plain}

\newpage
\tableofcontents

\end{document}